\documentclass[final,3p]{elsarticle}
\usepackage{amsfonts}
\usepackage{graphicx}
\usepackage{amsmath}
\usepackage{amssymb}
\usepackage{amsthm}
\usepackage{mathrsfs}
\usepackage{color}
\usepackage{makecell,multirow,diagbox}
\usepackage{amsmath}
\usepackage{amsfonts}
\usepackage{amsthm}
\usepackage{amsmath,amsxtra,latexsym,amsthm, amssymb, amscd, pb-diagram}
\usepackage{subfigure}

\newtheorem{dn}{Definition}[section]
\newtheorem{bdt}{Inequality}[section]
\newtheorem{dl}{Theorem}[section]
\newtheorem{md}{Proposition}[section]
\newtheorem{bd}{Lemma}[section]
\newtheorem{hq}{Corollary}[section]
\newtheorem{nx}{Remark}[section]
\newtheorem{vd}{Example}[section]
\newtheorem{tc}{Property}[section]
\newtheorem{cm}{Proof}[section]
\newcommand{\R}{\mathbb{R}}

\newcommand{\Z}{\mathbb{Z}}

\newcommand{\e}{\varepsilon}
\newcommand{\ity}{\infty}
\newcommand{\f}{\frac}

\newcommand{\bbd}{\begin{bd}}
\newcommand{\ebd}{\end{bd}}
\newcommand{\bbdt}{\begin{bdt}}
\newcommand{\ebdt}{\end{bdt}}
\newcommand{\bdn}{\begin{dn}}
\newcommand{\edn}{\end{dn}}
\newcommand{\bhq}{\begin{hq}}
\newcommand{\ehq}{\end{hq}}
\newcommand{\bdl}{\begin{dl}}
\newcommand{\edl}{\end{dl}}
\newcommand{\bnx}{\begin{nx}}
\newcommand{\enx}{\end{nx}}
\newcommand{\bmd}{\begin{md}}
\newcommand{\emd}{\end{md}}
\newcommand{\bcm}{\begin{cm}}
\newcommand{\ecm}{\end{cm}}
\newcommand{\bvd}{\begin{vd}}
\newcommand{\evd}{\end{vd}}
\newcommand{\btc}{\begin{tc}}
\newcommand{\etc}{\end{tc}}

\usepackage{epsfig}
\usepackage{enumitem}

\journal{Elsevier}

\begin{document}

\begin{frontmatter}




\author{Tuan Anh Dao\fnref{label1,label2}}
\ead{daotuananh.fami@gmail.com}
\author{Michael Reissig\corref{cor2}\fnref{label2}}
\ead{reissig@math.tu-freiberg.de}


\cortext[cor2]{Corresponding author.}


\title{$L^1$ estimates for oscillating integrals and their applications to semi-linear models with $\sigma$-evolution like structural damping}


\address[label1]{School of Applied Mathematics and Informatics, Hanoi University of Science and Technology, No.1 Dai Co Viet road, Hanoi, Vietnam}
\address[label2]{Faculty for Mathematics and Computer Science, TU Bergakademie Freiberg, Pr\"{u}ferstr. 9, 09596, Freiberg, Germany}


\begin{abstract}
The present paper is a continuation of our recent paper \cite{DaoReissig}. We will consider the following Cauchy problem for semi-linear structurally damped $\sigma$-evolution models:
\begin{equation*}
u_{tt}+ (-\Delta)^\sigma u+ \mu (-\Delta)^\delta u_t = f(u,u_t),\,\,\, u(0,x)= u_0(x),\,\,\,k u_t(0,x)=u_1(x)
\end{equation*}
with $\sigma \ge 1$, $\mu>0$ and $\delta \in (\frac{\sigma}{2},\sigma]$. Our aim is to study two main models including $\sigma$-evolution models with structural damping $\delta \in (\frac{\sigma}{2},\sigma)$ and those with visco-elastic damping $\delta=\sigma$. Here the function $f(u,u_t)$ stands for power nonlinearities $|u|^{p}$ and $|u_t|^{p}$ with a given number $p>1$. We are interested in investigating the global (in time) existence of small data Sobolev solutions to the above semi-linear model from suitable function spaces basing on $L^q$ space by assuming additional $L^{m}$ regularity for the initial data, with $q\in (1,\infty)$ and $m\in [1,q)$.
\end{abstract}


\end{frontmatter}


\section{Introduction}
\label{Introduction}
In the present paper, we study the following two Cauchy problems:
\begin{equation}
 u_{tt}+ (-\Delta)^\sigma u+ \mu (-\Delta)^{\delta} u_t=|u|^p,\,\,\, u(0,x)= u_0(x),\,\,\, u_t(0,x)=u_1(x), \label{pt6.1}
\end{equation}
and
\begin{equation}
 u_{tt}+ (-\Delta)^\sigma u+ \mu (-\Delta)^{\delta} u_t=|u_t|^p,\,\,\, u(0,x)= u_0(x),\,\,\, u_t(0,x)=u_1(x) \label{pt6.2}
\end{equation}
with $\sigma \ge 1$, $\mu>0$ and $\delta \in (\frac{\sigma}{2},\sigma]$. The corresponding linear model with vanishing right-hand side is
\begin{equation}
 u_{tt}+ (-\Delta)^\sigma u+ \mu (-\Delta)^{\delta} u_t=0,\,\,\, u(0,x)= u_0(x),\,\,\, u_t(0,x)=u_1(x). \label{pt6.3}
\end{equation}
A lot of papers (see, for example, \cite{DabbiccoEbert2014,DabbiccoReissig,NarazakiReissig}) focused on studying the special case $\sigma=1$ in (\ref{pt6.3}) with $\delta \in (0,1]$, that is, the model
\begin{equation}
 u_{tt}-\Delta u+ \mu (-\Delta)^{\delta} u_t=0,\,\,\, u(0,x)= u_0(x),\,\,\, u_t(0,x)=u_1(x). \label{pt6.4}
\end{equation}
More in detail, in the case $\delta \in (0,1)$ in \cite{NarazakiReissig} the authors  studied $L^1$ estimates for oscillating integrals to conclude $L^p- L^q$ estimates not necessarily on the conjugate line for solutions to (\ref{pt6.4}). In the case of semi-linear structurally damped wave models (\ref{pt6.1}) with $\sigma=1$ and $\delta \in (0,1]$ (see \cite{DabbiccoReissig}), the authors proved the global (in time) existence of small data solutions in low space dimensions by using classical energy estimates. In particular, they proposed to distinguish between ``parabolic like models" $\delta \in (0,\frac{1}{2})$ (see also \cite{MitidieriPohozaev}) and ``hyperbolic like models" $\delta \in (\frac{1}{2},1)$ (see also \cite{GalaktionovMitidieri}) from the point of decay estimates. Moreover, in \cite{DabbiccoEbert2014} some global (in time) existence results for small data solutions were presented for ``parabolic like models'' related to (\ref{pt6.1}) with $\sigma=1$ and $\delta \in (0,\frac{1}{2})$. More general, if we interested in studying (\ref{pt6.1}) and (\ref{pt6.2}) with $\delta \in [0,\frac{\sigma}{2}]$, then we want to mention the paper \cite{DabbiccoEbert} as another approach to obtain sharp $L^{p}- L^{q}$ estimates with $1<p \le q<\ity$ for the solutions to the linear model (\ref{pt6.3}). In detail, here the authors found an explicit way to get these estimates by using the Mikhlin-H\"{o}rmander multiplier theorem (see, for instance, \cite{Miyachi,RunSic}) for kernels localized at high frequencies. Then, there appeared some $L^{q}$ estimates for the solutions and some of their derivatives, with $q \in (1,\ity)$, to prove the global (in time) existence of small data solutions to the semi-linear models (\ref{pt6.1}) and (\ref{pt6.2}). In order to look for these results two different strategies were used due to the lack of $L^{1}- L^{1}$ estimates for solutions to (\ref{pt6.3}). They took account of additional $L^{1}\cap L^{\ity}$ regularity and additional $L^{\eta}\cap L^{\bar{q}}$ regularity for any small $\eta$ and large $\bar{q}$, respectively, in the first case with $\delta= \frac{\sigma}{2}$ and in the second case with $\delta \in (0,\frac{\sigma}{2})$. Recently, in \cite{DuongKainaneReissig} the use of $L^{2}-L^{2}$ estimates for solutions to (\ref{pt6.3}) by assuming additional $L^{1}$ regularity for the data was investigated to study semi-linear $\sigma$-evolution models (\ref{pt6.1}) and (\ref{pt6.2}) with $\delta= \frac{\sigma}{2}$. \medskip

Moreover, another interesting model related to (\ref{pt6.4}) is that with visco-elastic damping $\delta=1$ (or strong damping, see also \cite{Ikehata,IkehataTodorova}). It was considered in detail in \cite{Shibata}. The author obtained a potential decay estimate for solutions localized to low frequencies, whereas the high frequency part decays exponentially under the requirement of a suitable regularity for the data by application of the Marcinkiewicz theorem (see, for example, \cite{Marcinkiewicz, Weisz}) to related Fourier multipliers. The case of semi-linear visco-elastic damped wave models (\ref{pt6.1}) and (\ref{pt6.2}) with $\sigma=\delta=1$ was studied in several recent papers such as \cite{DabbiccoReissig} and \cite{Pizichillo}. In \cite{ReissigEbert} the authors mentioned some different interesting models related to (\ref{pt6.4}), namely those with $\sigma=\delta=2$, well-known as the visco-elastic damped plate models. Some decay estimates of the energy and qualitative properties of solutions as well were studied.\medskip

The present paper is a continuation of our recent paper \cite{DaoReissig}, in which the global (in time) existence of small data Sobolev solutions by mixing of additional $L^{m}$ regularity for the data on the basis of $L^{q}- L^{q}$ estimates, with $1 \le m< q< \infty$, is proved  to the semi-linear models (\ref{pt6.1}) and (\ref{pt6.2}). Here we remark that the properties of the solutions to (\ref{pt6.1}) and (\ref{pt6.2}) change completely from $(0,\frac{\sigma}{2})$ to $(\frac{\sigma}{2},\sigma]$. In particular, we want to distinguish between ``parabolic like models" ($\delta \in [0,\frac{\sigma}{2})$) (see \cite{DaoReissig}) and ``$\sigma$-evolution like models" ($\delta \in (\frac{\sigma}{2},\sigma]$) according to expected decay estimates. To do this, the first step of the present paper is to develop some $L^{1}$ estimates relying on several techniques from \cite{NarazakiReissig} for oscillating integrals in the presentation of solutions to (\ref{pt6.3}) by using theory of modified Bessel functions. It is also reasonable to apply Fa\`{a} di Bruno's formula (see, for instance, \cite{Bui,FrancescoBruno}) since the connection to Fourier multipliers appearing for wave models used in \cite{NarazakiReissig} fails to $\sigma$-evolution models for $\sigma >1$ (infinite speed of propagation of perturbations in the latter case). Then, we derive $L^{p}- L^{q}$ estimates not necessarily on the conjugate line for the solutions to (\ref{pt6.3}), with $1 \le p \le q \le \ity$, in the case of structural damping $\delta \in (\frac{\sigma}{2},\sigma)$. In the second step of this paper, we obtain $L^{q}- L^{q}$ estimates, with $1 \le q \le \infty$, for the solutions to (\ref{pt6.3}) by assuming suitable regularity for the data and applying the Mikhlin-H\"{o}rmander multiplier theorem for high frequencies in the remaining case of visco-elastic damping $\delta=\sigma$, for any $\sigma \ge 1$. Finally, having $L^{q}- L^{q}$ estimates by assuming additional $L^{m}$ regularity for the data and some developed tools from Harmonic Analysis in \cite{Palmierithesis} (see also \cite{DabbiccoReissig, Ozawa, Bui}) play a fundamental role to prove our global (in time) existence results. \medskip

The organization of this paper is as follows: \\
In Section \ref{Sec2}, we state the main results for the global (in time) existence of small data Sobolev solutions to (\ref{pt6.1}) and (\ref{pt6.2}). We present estimates for the solutions to (\ref{pt6.3}) in Section \ref{Sec3}. In particular, we provide estimates for solutions in the case of structural damping $\delta \in (\frac{\sigma}{2},\sigma)$ in Section \ref{Sec3.1} including the proof of $L^{1}$ estimates, $L^{\infty}$ estimates and $L^{r}$ estimates as well. Section \ref{Sec3.2} is devoted to derive  estimates of solutions in the case of visco-elastic damping $\delta=\sigma$. Then, in Section \ref{Sec3.3} we state $L^{q}- L^{q}$ estimates by assuming additional $L^{m}$ regularity for the data with $q\in (1,\infty)$ and $m\in [1,q)$. Then, we prove our global (in time) existence results to (\ref{pt6.1}) and (\ref{pt6.2}) in Section \ref{Sec4}. Finally, in Section \ref{Sec5} we state some concluding remarks and open problems.\medskip

Throughout the present paper, we use the following notations.\medskip

\noindent\textbf{Notation 1.} We write $f\lesssim g$ when there exists a constant $C>0$ such that $f\le Cg$, and $f \approx g$ when $g\lesssim f\lesssim g$.\medskip

\noindent\textbf{Notation 2.} We denote $[s]^+:= \max\{s,0\}$ as the positive part of $s \in \R$, and $\lceil s \rceil:= \min \big\{k \in \Z \,\, : \,\, k\ge s \big\}$.\medskip

\noindent\textbf{Notation 3.} The spaces $H^{a,q}$ and $\dot{H}^{a,q}$, with $a \ge 0$ and $q> 1$, denote Bessel and Riesz potential spaces based on $L^q$. As usual, $\big<D\big>^{a}$ and $|D|^{a}$ stand for the pseudo-differential operators with symbols $\big<\xi\big>^{a}$ and $|\xi|^{a}$, respectively.\medskip

\noindent\textbf{Notation 4.} We introduce the space $\mathcal{A}^{s}_{m,q}:= \big(L^m \cap H^{s,q}\big) \times \big(L^m \cap H^{[s-2\delta]^+,q}\big)$ with the norm
$$\|(u_0,u_1)\|_{\mathcal{A}^{s}_{m,q}}:=\|u_0\|_{L^m}+ \|u_0\|_{H^{s,q}}+ \|u_1\|_{L^m}+ \|u_1\|_{H^{[s-2\delta]^+,q}}, \text{ for } s\ge 0. $$

\noindent\textbf{Notation 5.} We fix the constants $\kappa_1:= 1+(1+[\frac{n}{2}])(1-\frac{\sigma}{2\delta})\big(1+\frac{1}{q}-\frac{1}{m}\big)$ and $\kappa_2:= (2+[\frac{n}{2}])(1-\frac{\sigma}{2\delta})\big(1+\frac{1}{q}-\frac{1}{m}\big)$.

\section{Main results} \label{Sec2}
In the first case, we obtain solutions to (\ref{pt6.1}) from energy space on the base of the space $L^q$.
\bdl \label{dl2.1}
Let $q \in (1,\ity)$ be a fixed constant and $m\in [1,q)$ in (\ref{pt6.1}). We assume the condition
\begin{equation}
p> 1+\f{\max\{2m\delta(1+\kappa_1),\, n-\frac{m}{q}n+2m\delta\}}{n-2m\delta\kappa_1}. \label{exponent2A}
\end{equation}
Moreover, we suppose the following conditions:
\begin{equation}
p \in \Big[\f{q}{m}, \ity \Big) \text{ if } n \le 2q\delta, \text{ or }p \in \Big[\f{q}{m}, \f{n}{n-2q\delta}\Big] \text{ if } n \in \Big(2q\delta, \f{2q^2 \delta}{q-m}\Big]. \label{GN2A}
\end{equation}
Then, there exists a constant $\e>0$ such that for any small data
$$(u_0,u_1) \in \mathcal{A}^{2\delta}_{m,q} \text{ satisfying the assumption } \|(u_0,u_1)\|_{\mathcal{A}^{2\delta}_{m,q}} \le \e, $$
we have a uniquely determined global (in time) small data energy solution (on the base of $L^q$) \[ u\in C([0,\ity),H^{2\delta,q})\cap C^1([0,\ity),L^q)\]
to (\ref{pt6.1}). The following estimates hold:
\begin{align}
\big\|u(t,\cdot)\big\|_{L^q}& \lesssim (1+t)^{1+(1+[\frac{n}{2}])(1-\frac{\sigma}{2\delta})\frac{1}{r}-\frac{n}{2\delta}(1-\frac{1}{r})} \|(u_0,u_1)\|_{\mathcal{A}^{2\delta}_{m,q}}, \label{decayrate2A1}\\
\big\||D|^{\sigma}u(t,\cdot)\big\|_{L^q}& \lesssim (1+t)^{1+(1+[\frac{n}{2}])(1-\frac{\sigma}{2\delta})\frac{1}{r} -\frac{n}{2\delta}(1-\frac{1}{r})- \frac{\sigma}{2\delta}}\|(u_0,u_1)\|_{\mathcal{A}^{2\delta}_{m,q}}, \label{decayrate2A2}\\
\big\|u_t(t,\cdot)\big\|_{L^q}& \lesssim (1+t)^{(2+[\frac{n}{2}])(1-\frac{\sigma}{2\delta})\frac{1}{r} -\frac{n}{2\delta}(1-\frac{1}{r})} \|(u_0,u_1)\|_{\mathcal{A}^{2\delta}_{m,q}}, \label{decayrate2A3}\\
\big\||D|^{2\delta}u(t,\cdot)\big\|_{L^q}& \lesssim (1+t)^{(1+[\frac{n}{2}])(1-\frac{\sigma}{2\delta})\frac{1}{r} -\frac{n}{2\delta}(1-\frac{1}{r})}\|(u_0,u_1)\|_{\mathcal{A}^{2\delta}_{m,q}}, \label{decayrate2A4}
\end{align}
where $1+\frac{1}{q}=\frac{1}{r}+\frac{1}{m}$.
\edl

In the second case, we obtain Sobolev solutions to (\ref{pt6.1}).

\bdl \label{dl2.2}
Let $q \in (1,\ity)$ be a fixed constant, $m\in [1,q)$ in (\ref{pt6.1}) and $0 < s < 2\delta$. We assume the condition
\begin{equation}
p> 1+\f{\max\{2m\delta(1+\kappa_1),\, n-\frac{m}{q}n+ms\}}{n-2m\delta\kappa_1}. \label{exponent3A}
\end{equation}
Moreover, we suppose the following conditions:
\begin{equation}
p \in \Big[\f{q}{m}, \ity \Big) \text{ if } n \le qs, \text{ or }p \in \Big[\f{q}{m}, \f{n}{n-qs}\Big] \text{ if } n \in \Big(qs, \f{q^2 s}{q-m}\Big]. \label{GN3A}
\end{equation}
Then, there exists a constant $\e>0$ such that for any small data
$$(u_0,u_1) \in \mathcal{A}^{s}_{m,q} \text{ satisfying the assumption } \|(u_0,u_1)\|_{\mathcal{A}^{s}_{m,q}}\le \e, $$
we have a uniquely determined global (in time) small data Sobolev solution \[ u\in C([0,\ity),H^{s,q})\]
to (\ref{pt6.1}). The following estimates hold:
\begin{align}
\big\|u(t,\cdot)\big\|_{L^q}& \lesssim (1+t)^{1+(1+[\frac{n}{2}])(1-\frac{\sigma}{2\delta})\frac{1}{r}-\frac{n}{2\delta}(1-\frac{1}{r})} \|(u_0,u_1)\|_{\mathcal{A}^{s}_{m,q}}, \label{decayrate3A1}\\
\big\||D|^{s}u(t,\cdot)\big\|_{L^q}& \lesssim (1+t)^{1+(1+[\frac{n}{2}])(1-\frac{\sigma}{2\delta})\frac{1}{r} -\frac{n}{2\delta}(1-\frac{1}{r})-\frac{s}{2\delta}} \|(u_0,u_1)\|_{\mathcal{A}^{s}_{m,q}}, \label{decayrate3A2}
\end{align}
where $1+\frac{1}{q}=\frac{1}{r}+\frac{1}{m}$.
\edl

\begin{nx}
\fontshape{n}
\selectfont
We want to underline that due to the flexibility of parameter $q \in (1,\ity)$, we really get a result for arbitrarily small positive $s$ in Theorem \ref{dl2.2}. In particular, if we take any small positve $s=\e$, then we also choose for example a sufficiently large $q=\frac{1}{\e^2}$ in order to guarantee the existence of both the space dimension $n$ and the exponent $p$ satisfying the required conditions in Theorem \ref{dl2.2}.
\end{nx}

In the third case, we obtain solutions to (\ref{pt6.1}) belonging to the energy space (on the base of $L^q$) with a suitable higher regularity.

\bdl \label{dl2.3}
Let $q \in (1,\ity)$ be a fixed constant, $m\in [1,q)$ in (\ref{pt6.1}) and $2\delta<s \le 2\delta+\frac{n}{q}$. We assume the exponent $p>1+ \lceil s- 2\delta \rceil$ satisfying the condition
\begin{equation}
p> 1+\f{\max\{2m\delta(1+\kappa_1),\, n-\frac{m}{q}n+ms\}}{n-2m\delta\kappa_1}. \label{exponent4A}
\end{equation}
Moreover, we suppose the following conditions:
\begin{equation}
p \in \Big[\frac{q}{m}, \ity \Big) \text{ if } n \le qs, \text{ or }p \in \Big[\frac{q}{m}, 1+\f{2q\delta}{n-qs}\Big] \text{ if }  \Big(qs, qs+\f{2mq\delta}{q-m}\Big]. \label{GN4A}
\end{equation}
Then, there exists a constant $\e>0$ such that for any small data
$$(u_0,u_1) \in \mathcal{A}^{s}_{m,q}\text{ satisfying the assumption }\|(u_0,u_1)\|_{\mathcal{A}^{s}_{m,q}}\le \e, $$
we have a uniquely determined global (in time) small data energy solution \[ u\in C([0,\ity),H^{s,q})\cap C^1([0,\ity),H^{s-2\delta,q})\]
to (\ref{pt6.1}). The following estimates hold:
\begin{align}
\|u(t,\cdot)\|_{L^q}& \lesssim (1+t)^{1+(1+[\frac{n}{2}])(1-\frac{\sigma}{2\delta})\frac{1}{r}-\frac{n}{2\delta}(1-\frac{1}{r})} \|(u_0,u_1)\|_{\mathcal{A}^{s}_{m,q}}, \label{decayrate4A1}\\
\big\||D|^s u(t,\cdot)\big\|_{L^q}& \lesssim (1+t)^{1+(1+[\frac{n}{2}])(1-\frac{\sigma}{2\delta})\frac{1}{r} -\frac{n}{2\delta}(1-\frac{1}{r})-\frac{s}{2\delta}} \|(u_0,u_1)\|_{\mathcal{A}^{s}_{m,q}}, \label{decayrate4A2}\\
\|u_t(t,\cdot)\|_{L^q}& \lesssim (1+t)^{(2+[\frac{n}{2}])(1-\frac{\sigma}{2\delta})\frac{1}{r} -\frac{n}{2\delta}(1-\frac{1}{r})} \|(u_0,u_1)\|_{\mathcal{A}^{s}_{m,q}}, \label{decayrate4A3}\\
\big\||D|^{s-2\delta}u_t(t,\cdot)\big\|_{L^q}& \lesssim (1+t)^{1+(2+[\frac{n}{2}])(1-\frac{\sigma}{2\delta})\frac{1}{r} -\frac{n}{2\delta}(1-\frac{1}{r}) -\frac{s}{2\delta}} \|(u_0,u_1)\|_{\mathcal{A}^{s}_{m,q}}, \label{decayrate4A4}
\end{align}
where $1+\frac{1}{q}=\frac{1}{r}+\frac{1}{m}$.
\edl

\begin{nx}
\fontshape{n}
\selectfont
Let us explain the conditions for $p$ and $n$ in Theorems \ref{dl2.1} to \ref{dl2.3}. The conditions (\ref{exponent2A}), (\ref{exponent3A}) and (\ref{exponent4A}) imply the same decay estimates for the solutions to (\ref{pt6.1}) as for the solutions to the corresponding linear Cauchy problem (\ref{pt6.3}). Hence, we can say that the nonlinearity is interpreted as a small perturbation. The other conditions (\ref{GN2A}) and (\ref{GN3A}) come into play after we apply fractional Gargliardo-Nirenberg inequality. In addition, the upper bound for $n$ arises from the corresponding set of admissible parameters $p$ to guarantee that this range is non-empty. Employing fractional chain rule leads to the condition $p>1+ \lceil s- 2\delta \rceil$ in Theorem \ref{dl2.3}. Eventually, the condition (\ref{GN4A}) appears as an interplay between fractional Gargliardo-Nirenberg inequality and fractional chain rule.
\end{nx}
Finally, we obtain high regular solutions to (\ref{pt6.1}) by using the fractional Sobolev embedding.
\bdl \label{dl2.4}
Let $s >2\delta+\frac{n}{q}$. Let $q \in (1,\ity)$ be a fixed constant and $m\in [1,q)$ in (\ref{pt6.1}). We assume that the exponent $p> 1+ s- 2\delta$ satisfies the condition
\begin{equation}
p> 1+\f{\max\{2m\delta(1+\kappa_1),\,\,\, n-\frac{m}{q}n+ms\}}{n-2m\delta\kappa_1}. \label{exponent5A}
\end{equation}
Moreover, we suppose the following conditions:
\begin{equation}
p \in \Big[\frac{q}{m}, \ity \Big) \text{ and } n> 2m\delta\kappa_1. \label{GN5A}
\end{equation}
Then, there exists a constant $\e>0$ such that for any small data
$$(u_0,u_1) \in \mathcal{A}^{s}_{m,q}\text{ satisfying the assumption }\|(u_0,u_1)\|_{\mathcal{A}^{s}_{m,q}}\le \e, $$
we have a uniquely determined global (in time) small data energy solution \[ u\in C([0,\ity),H^{s,q})\cap C^1([0,\ity),H^{s-2\delta,q})\]
to (\ref{pt6.1}). Moreover, the estimates (\ref{decayrate4A1}) to (\ref{decayrate4A4}) hold.
\edl
Finally, we obtain large regular solutions to (\ref{pt6.2}) by using the fractional Sobolev embedding.
\bdl \label{dl2.5}
Let $s >2\delta+\frac{n}{q}$. Let $q \in (1,\ity)$ be a fixed constant and $m\in [1,q)$ in (\ref{pt6.2}). We assume that the exponent $p> 1+ s- 2\delta$ satisfies the condition
\begin{equation}
p> 1+\f{\max\{2m\delta(1+\kappa_2),\,\,\, n-\frac{m}{q}n+m(s-\sigma)\}}{n-2m\delta\kappa_2}. \label{exponent6A}
\end{equation}
Moreover, we suppose the following conditions:
\begin{equation}
p \in \Big[\frac{q}{m}, \ity \Big) \text{ and } n> 2m\delta\kappa_2. \label{GN6A}
\end{equation}
Then, there exists a constant $\e>0$ such that for any small data
$$(u_0,u_1) \in \mathcal{A}^{s}_{m,q}\text{ satisfying the assumption }\|(u_0,u_1)\|_{\mathcal{A}^{s}_{m,q}}\le \e, $$
we have a uniquely determined global (in time) small data energy solution \[ u\in C([0,\ity),H^{s,q})\cap C^1([0,\ity),H^{s-2\delta,q})\]
to (\ref{pt6.2}). Moreover, the estimates (\ref{decayrate4A1}) to (\ref{decayrate4A4}) hold.
\edl
\begin{nx}
\fontshape{n}
\selectfont
Let us turn to interpret the conditions for $p$ and $n$ in Theorems \ref{dl2.4} and \ref{dl2.5}. Since we want to use fractional powers, the condition $p> 1+ s-2\delta$ is necessary. Moreover, the conditions (\ref{exponent5A}) and (\ref{exponent6A}) bring the same decay estimates for the solutions, respectively, to (\ref{pt6.1}) and (\ref{pt6.2}) as for solutions to the corresponding Cauchy problem (\ref{pt6.3}). Hence, the nonlinearity can be considered as a small perturbation. Finally, the remaining conditions (\ref{GN5A}) and (\ref{GN6A}) come from applying fractional Gargliardo-Nirenberg inequality and fractional Sobolev embedding.
\end{nx}
\begin{nx}
\fontshape{n}
\selectfont
In comparison with all the theorems in our previous paper \cite{DaoReissig}, we want to underline that the solutions from all the above theorems have no loss of regularity (see also \cite{DabbiccoEbert,Miyachi1980,Peral}) with respect to the initial data. Loss of regularity of the solutions appearing in \cite{DaoReissig} is due to the singular behavior of time-dependent coefficients in the estimates of solutions to the linear models localized to high frequencies as $t\longrightarrow +0$ with $\delta \in (0,\frac{\sigma}{2})$. Meanwhile, this phenomenon does not appear in the case  $\delta \in (\frac{\sigma}{2},\sigma]$ (see later, Proposition \ref{md4.3.3}).
\end{nx}
\begin{nx}
\fontshape{n}
\selectfont
Let us compare our results with some known results from \cite{DabbiccoReissig}. By choosing $\sigma=1$, $q=2$ and $m=1$ we see that, on the one hand, the admissible exponents $p$ in the cited paper are somehow better than those in Theorem \ref{dl2.1} for low space dimensions. On the other hand, we want to underline that Theorem \ref{dl2.1} allows some flexibility for both $p$ and $n$ because of the flexible choice of parameters $\sigma$, $\delta$, $q$ and $m$ (see also some of the examples below).
\end{nx}
\begin{vd}
\fontshape{n}
\selectfont
In the following examples, we choose $m=1$, $q=4$, $\sigma=\frac{9}{5}$, $\delta=1$ and $n=4$:
\begin{itemize}
\item Using Theorem \ref{dl2.1} we obtain $p \in [4,\ity)$.
\item If $s= \frac{3}{2}$, then using Theorem \ref{dl2.2} we obtain $p \in [4,\ity)$.
\item If $s= \frac{5}{2}$, then using Theorem \ref{dl2.3} we obtain $p \in [4,\ity)$.
\item If $s= \frac{7}{2}$, then using Theorem \ref{dl2.4} we obtain $p \in \big(\frac{167}{37},\ity\big)$.
\item If $s= 4$, then using Theorem \ref{dl2.5} we obtain $p \in [4,\ity)$.
\end{itemize}
\end{vd}
\begin{vd}
\fontshape{n}
\selectfont
In the following examples, we choose $m=1$, $q=4$, $\sigma=\delta=\frac{11}{10}$ and $n=5$:
\begin{itemize}
\item Using Theorem \ref{dl2.1} we obtain $p \in \big(\frac{317}{79},\ity\big)$.
\item If $s= 2$, then using Theorem \ref{dl2.2} we obtain $p \in [4,\ity)$.
\item If $s= \frac{5}{2}$, then using Theorem \ref{dl2.3} we obtain $p \in \big(\frac{329}{79},\ity\big)$.
\item If $s= \frac{7}{2}$, then using Theorem \ref{dl2.4} we obtain $p \in \big(\frac{369}{79},\ity\big)$.
\item If $s= 4$, then using Theorem \ref{dl2.5} we obtain $p \in [4,\ity)$.
\end{itemize}
\end{vd}

\section{Decay estimates for solutions to linear Cauchy problems} \label{Sec3}
The goal of this section is to obtain decay estimates for the solution and some its derivatives to (\ref{pt6.3}). These estimates play an essential role to prove the global (in time) existence results to (\ref{pt6.1}) and (\ref{pt6.2}) in the next section. First, using partial Fourier transformation to (\ref{pt6.3}) leads to the following Cauchy problem for $v(t,\xi):=F_{x\rightarrow \xi}\big(u(t,x)\big)$, $v_0(\xi):=F_{x\rightarrow \xi}\big(u_0(x)\big)$ and $v_1(\xi):=F_{x\rightarrow \xi}\big(u_1(x)\big)$
\begin{equation}
v_{tt}+ \mu |\xi|^{2\delta} v_t+ |\xi|^{2\sigma} v=0,\,\, v(0,\xi)= v_0(\xi),\,\, v_t(0,\xi)= v_1(\xi). \label{pt3.1}
\end{equation}
Without loss of generality we can choose $\mu=1$ in (\ref{pt3.1}). The characteristic roots are
$$ \lambda_{1,2}=\lambda_{1,2}(\xi)= \f{1}{2}\Big(-|\xi|^{2\delta}\pm \sqrt{|\xi|^{4\delta}-4|\xi|^{2\sigma}}\Big). $$
The solution to (\ref{pt3.1}) can be written as follows:
\begin{equation}
v(t,\xi)= \frac{\lambda_1 e^{\lambda_2 t}-\lambda_2 e^{\lambda_1 t}}{\lambda_1- \lambda_2}v_0(\xi)+ \frac{e^{\lambda_1 t}-e^{\lambda_2 t}}{\lambda_1- \lambda_2}v_1(\xi)=: \hat{K_0}(t,\xi)v_0(\xi)+\hat{K_1}(t,\xi)v_1(\xi). \label{pt3.2}
\end{equation}
Here we assume $\lambda_{1}\neq \lambda_{2}$. Taking account of the cases of small and large frequencies separately, we get
\begin{align}
&1.\,\, \lambda_{1,2} \sim -|\xi|^{2\delta}\pm i|\xi|^\sigma,\,\, \lambda_1-\lambda_2 \sim i|\xi|^\sigma \text{ for small } |\xi|, \label{pt3.3}\\
&2.\,\, \lambda_1\sim -|\xi|^{2(\sigma- \delta)},\,\, \lambda_2\sim -|\xi|^{2\delta},\,\, \lambda_1-\lambda_2 \sim |\xi|^{2\delta} \text{ for large } |\xi|. \label{pt3.4}
\end{align}
We now decompose the solution to (\ref{pt6.3}) into two parts localized separately at low and high frequencies, that is,
$$ u(t,x)= u_\chi(t,x)+ u_{1-\chi}(t,x), $$
where
$$u_\chi(t,x)=F^{-1}\big(\chi(|\xi|)v(t,\xi)\big) \text{ and } u_{1-\chi}(t,x)=F^{-1}\Big(\big(1-\chi(|\xi|)\big)v(t,\xi)\Big), $$
with a smooth cut-off function $\chi=\chi(|\xi|)$ equal to $1$ for small $|\xi|$ and vanishing for large $|\xi|$.
\subsection{Estimates for oscillating integrals in the case of structural damping: $\delta \in (\frac{\sigma}{2},\sigma)$} \label{Sec3.1}
\subsubsection{$L^1$ estimates for small frequencies}  \label{Sec3.1.1}

\bmd \label{md4.1.2}
The following estimates hold in $\R^n$ for any $n \ge 1$:
\begin{align*}
&\big\| F^{-1}\big(|\xi|^a \hat{K_0}\chi(|\xi|)\big)(t,\cdot)\big\|_{L^1}\lesssim \begin{cases}
1 \text{ for } t\in (0,1], &\\
t^{(2+[\frac{n}{2}])(1-\frac{\sigma}{2\delta})-\frac{a}{2\delta}} \text{ for } t\in[1,\ity), &
\end{cases}\\
&\big\| F^{-1}\big(|\xi|^a \hat{K_1}\chi(|\xi|)\big)(t,\cdot)\big\|_{L^1}\lesssim \begin{cases}
t \text{ for } t\in (0,1], &\\
t^{1+(1+[\frac{n}{2}])(1-\frac{\sigma}{2\delta})-\frac{a}{2\delta}} \text{ for } t\in[1,\ity), &
\end{cases}
\end{align*}
for any non-negative number $a$.
\emd
To derive the desired estimates for the norm of the Fourier multipliers localized to small frequencies, we write
\begin{equation}
|\xi|^a \hat{K_0}(t,\xi) =e^{-\frac{1}{2}|\xi|^{2\delta}t}|\xi|^a \cos \Big( |\xi|^\sigma \sqrt{1-\frac{1}{4}|\xi|^{4\delta- 2\sigma}} t \Big)+ e^{-\frac{1}{2} |\xi|^{2\delta}t}|\xi|^{a+2\delta}\frac{\sin \big( |\xi|^\sigma \sqrt{1-\frac{1}{4}|\xi|^{4\delta- 2\sigma}}t \big)}{2|\xi|^\sigma \sqrt{1-\frac{1}{4}|\xi|^{4\delta- 2\sigma}}}, \label{p3.2.1}
\end{equation}
and
\begin{equation}
|\xi|^a \hat{K_1}(t,\xi)=e^{-\frac{1}{2} |\xi|^{2\delta}t}|\xi|^a \frac{\sin \big(|\xi|^\sigma \sqrt{1-\frac{1}{4}|\xi|^{4\delta- 2\sigma}}t \big)}{|\xi|^\sigma \sqrt{1-\frac{1}{4}|\xi|^{4\delta- 2\sigma}}}. \label{p3.2.2}
\end{equation}
For this reason, we will split our proof into two steps. In the first step we derive $L^1$ estimates for the following oscillating integrals:
$$ F^{-1}\Big(e^{-c_1 |\xi|^{2\delta}t}|\xi|^{2\beta}\frac{\sin (c_2 |\xi|^\sigma t)}{|\xi|^\sigma}\chi(|\xi|)\Big)(t,\cdot), $$
and
$$ F^{-1}\Big(e^{-c_1 |\xi|^{2\delta}t}|\xi|^{2\beta}\cos (c_2 |\xi|^\sigma t) \chi(|\xi|)\Big)(t,\cdot), $$
where $\beta \ge 0$, $c_1$ is a positive constant and $c_2 \ne 0$ is a real constant. Then, in the second step we estimate the following more structured oscillating integrals:
$$ F^{-1}\Big(e^{-c_1 |\xi|^{2\delta}t}|\xi|^{2\beta}\frac{\sin \big(c_2 |\xi|^\sigma f(|\xi|) t\big)}{|\xi|^\sigma f(|\xi|)}\chi(|\xi|)\Big)(t,\cdot), $$
and
$$ F^{-1}\Big(e^{-c_1 |\xi|^{2\delta}t}|\xi|^{2\beta}\cos \big(c_2 |\xi|^\sigma f(|\xi|) t\big) \chi(|\xi|)\Big)(t,\cdot), $$
where $$f(|\xi|)=\sqrt{1-\frac{1}{4}|\xi|^{4\delta- 2\sigma}}. $$

\bbd \label{bd4.1.1}
The following estimate holds in $\R^n$ for any $n \ge 1$:
$$\Big\| F^{-1}\Big(e^{-c_1 |\xi|^{2\delta}t}|\xi|^{2\beta}\frac{\sin (c_2 |\xi|^\sigma t)}{|\xi|^\sigma}\chi(|\xi|)\Big)(t,\cdot)\Big\|_{L^1} \lesssim \begin{cases}
t \text{ for } t\in (0,1], &\\
t^{(2+[\frac{n}{2}])(1-\frac{\sigma}{2\delta})+\frac{\sigma-2\beta}{2\delta}} \text{ for } t\in[1,\ity), &
\end{cases}$$
with $\beta \ge 0$. Here $c_1$ is a positive and $c_2 \ne 0$ is a real constant.
\ebd
\begin{proof} We follow ideas from the proofs to Proposition $4$ in \cite{NarazakiReissig} and Lemma $3.1$ in \cite{DaoReissig}. Many steps in our proof are similar to the proofs of these results. Hence, it is reasonable to present only the steps which are different. Let us divide the proof into two cases: $t\in (0,1]$ and $t\in [1,\ity)$. First, in order to treat the first case $t\in (0,1]$, we localize to small $|x|\le 1$. Then, we derive immediately for small values of $|\xi|$ the estimate
\begin{equation}
\Big\|F^{-1}\Big(e^{-c_1 |\xi|^{2\delta}t}|\xi|^{2\beta-\sigma} \sin (c_2 |\xi|^\sigma t) \chi(|\xi|)\Big)(t,\cdot)\Big\|_{L^1(|x| \le 1)} \lesssim t. \label{l4.1.1}
\end{equation}
For this reason we assume now $|x| \ge 1$. We introduce the function
$$I=I(t,x):=F^{-1}\Big(e^{-c_1 |\xi|^{2\delta}t}|\xi|^{2\beta-\sigma} \sin (c_2 |\xi|^\sigma t) \chi(|\xi|)\Big)(t,x). $$
Because the functions in the parenthesis are radial symmetric with respect to $\xi$, the inverse Fourier transform is radial symmetric with respect to $x$, too. Applying modified Bessel functions leads to
\begin{equation}
I(t,x)= c\int_0^\ity e^{-c_1 r^{2\delta}t}r^{2\beta-\sigma} \sin (c_2 r^\sigma t) \chi(r) r^{n-1} \tilde{J}_{\frac{n}{2}-1}(r|x|)\, dr. \label{l4.1.2}
\end{equation}
Let us consider odd spatial dimension $n=2m+1, m \ge 1$. We introduce the vector field $Xf(r):= \frac{d}{dr}\big(\frac{1}{r}f(r)\big)$ as in the proof of Proposition $4$ in \cite{NarazakiReissig}. Then carrying out $m+1$ steps of partial integration we have
\begin{equation}
I(t,x)=-\frac{c}{|x|^{n}}\int_0^\ity \partial_r \big(X^m \big( e^{-c_1 r^{2\delta}t} \sin (c_2 r^\sigma t)\chi(r) r^{2\beta-\sigma+2m}\big)\big) \sin(r|x|)\,dr. \label{l4.1.3}
\end{equation}
A standard calculation leads to the following presentation of the right-hand side of (\ref{l4.1.3}):
\begin{align*}
I(t,x)&= \sum_{j=0}^m \sum_{k=0}^{j+1}\frac{c_{jk}}{|x|^{n}}\int_0^\ity \partial_r^{j+1-k} e^{-c_1 r^{2\delta}t} \partial_r^{k} \big(\sin (c_2 r^\sigma t)\chi(r)\big) r^{2\beta-\sigma+j} \sin(r|x|)\,dr\\
&\quad + \sum_{j=0}^m \sum_{k=0}^j \frac{c_{jk}}{|x|^{n}}\int_0^\ity \partial_r^{j-k} e^{-c_1 r^{2\delta}t}\partial_r^{k+1} \big(\sin (c_2 r^\sigma t)\chi(r)\big) r^{2\beta-\sigma+j} \sin(r|x|)\,dr\\
&\quad + \sum_{j=1}^m \sum_{k=0}^j \frac{c_{jk}}{|x|^{n}}\int_0^\ity \partial_r^{j-k} e^{-c_1 r^{2\delta}t} \partial_r^{k} \big(\sin (c_2 r^\sigma t)\chi(r)\big) r^{2\beta-\sigma+j-1} \sin(r|x|)\,dr
\end{align*}
with some constants $c_{jk}$. Now, we control the integrals
\begin{equation}
I_{j,k}(t,x):= \int_0^\ity \partial_r^{j+1-k} e^{-c_1 r^{2\delta}t} \partial_r^{k} \big(\sin (c_2 r^\sigma t)\chi(r)\big) r^{2\beta-\sigma+j} \sin(r|x|)\,dr. \label{l4.1.4}
\end{equation}
Because of small values of $r$, we notice that the following estimates hold on the support of $\chi$ and on the support of its derivatives:
\begin{align*}
\big|\partial_r^l e^{-c_1 r^{2\delta}t}\big| &\lesssim
\begin{cases}
1 \text{ if } l=0, &\\
r^{2\delta-l} t \text{ if } l=1,\cdots,m,&
\end{cases}\\
\big|\partial_r^{l}\big(\sin (c_2 r^\sigma t)\chi(r)\big)\big| &\lesssim r^{\sigma-l}t \text{ for all } l=0,\cdots,m.
\end{align*}
As a result, we obtain for small $r$, $j=0,\cdots,m$ and $k=0,\cdots,j$
$$\big|\partial_r^{j+1-k} e^{-c_1 r^{2\delta}t} \partial_r^{k} \big(\sin (c_2 r^\sigma t)\chi(r)\big) r^{2\beta-\sigma+j}\big| \lesssim r^{2\delta+2\beta-1}t^2$$
on the support of $\chi$ and on the support of its derivatives. We divide the integral (\ref{l4.1.4}) into two parts to derive on the one hand
\begin{equation}
\Big|\int_0^{\frac{\pi}{2|x|}} \partial_r^{j+1-k} e^{-c_1 r^{2\delta}t} \partial_r^{k}\big(\sin (c_2 r^\sigma t)\chi(r)\big) r^{2\beta-\sigma+j} \sin(r|x|)\,dr \Big| \lesssim \frac{t^2}{|x|^{2\delta}}. \label{l4.1.5}
\end{equation}
On the other hand, we can carry out one more step of partial integration in estimating the remaining integral as follows:
\begin{align}
&\Big|\int_{\frac{\pi}{2|x|}}^\ity \partial_r^{j+1-k} e^{-c_1 r^{2\delta}t} \partial_r^{k}\big(\sin (c_2 r^\sigma t)\chi(r)\big) r^{2\beta-\sigma+j} \sin(r|x|)\,dr \Big| \nonumber \\
&\quad \lesssim \frac{1}{|x|}\Big|\partial_r^{j+1-k} e^{-c_1 r^{2\delta}t} \partial_r^{k}\big(\sin (c_2 r^\sigma t)\chi(r)\big) r^{2\beta-\sigma+j} \cos(r|x|) \Big|_{r=\frac{\pi}{2|x|}}^\ity \nonumber \\
&\qquad + \frac{1}{|x|}\int_{\frac{\pi}{2|x|}}^\ity \Big|\partial_r \Big( \partial_r^{j+1-k} e^{-c_1 r^{2\delta}t} \partial_r^{k}\big(\sin (c_2 r^\sigma t)\chi(r)\big) r^{2\beta-\sigma+j}\Big) \cos(r|x|)\Big|\, dr \lesssim \frac{t^2}{|x|}, \label{l4.1.6}
\end{align}
since $\delta+2\beta > \sigma \ge 1$. Here we also note that for all $j=0,\cdots,m$ and $k=0,\cdots,j$ we have
$$\Big|\partial_r \Big( \partial_r^{j+1-k} e^{-c_1 r^{2\delta}t} \partial_r^{k}\big(\sin (c_2 r^\sigma t)\chi(r)\big) r^{2\beta-\sigma+j}\Big)\Big| \lesssim r^{2\delta+2\beta-2}t^2. $$
Hence, from (\ref{l4.1.3}) to (\ref{l4.1.6}) we have produced terms $|x|^{-(n+2\delta)}$ and $|x|^{-(n+1)}$ which guarantee the $L^1$ property in $x$ to prove that for all $t \in (0,1]$ and $n=2m+1$ the following estimate holds:
\begin{equation}
\big\|F^{-1}\big(e^{-c_1 |\xi|^{2\delta}t}|\xi|^{2\beta-\sigma} \sin (c_2 |\xi|^\sigma t) \chi(|\xi|)\big)(t,\cdot)\big\|_{L^1(|x| \ge 1)} \lesssim t^2. \label{l4.1.7}
\end{equation}
Let us consider even spatial dimension $n=2m,\, m \ge 1$, in the first case $t \in (0,1]$. Then, applying the first rule of modified Bessel functions for $\mu=1$ and the fifth rule for $\mu=0$, and repeating the above calculations as we did to get (\ref{l4.1.7}) we can conclude the following estimate:
\begin{equation}
\big\|F^{-1}\big(e^{-c_1 |\xi|^{2\delta}t}|\xi|^{2\beta-\sigma} \sin (c_2 |\xi|^\sigma t) \chi(|\xi|)\big)(t,\cdot)\big\|_{L^1(|x| \ge 1)} \lesssim t^2, \text{ for } n=2m,\, m \ge 1. \label{l4.1.8}
\end{equation}
\noindent Let us turn to the second case $t \in [1,\ity)$. Then, by the change of variables $\xi=t^{-\frac{1}{2\delta}}\eta $ as we did in the proof of the case $t\in (0,1]$ to Lemma $3.1$ in \cite{DaoReissig} we will follow the steps of the proof of this lemma to conclude the following estimates:
\begin{equation}
\big\|F^{-1}\big(e^{-c_1 |\xi|^{2\delta}t}|\xi|^{2\beta-\sigma} \sin (c_2 |\xi|^\sigma t) \chi(|\xi|)\big)(t,\cdot)\big\|_{L^1(|x| \le 1)} \lesssim t^{1-\frac{\beta}{\delta}}, \label{l4.1.9}
\end{equation}
and
\begin{equation}
\big\|F^{-1}\big(e^{-c_1 |\xi|^{2\delta}t}|\xi|^{2\beta-\sigma} \sin (c_2 |\xi|^\sigma t) \chi(|\xi|)\big)(t,\cdot)\big\|_{L^1(|x| \ge 1)}\lesssim
\begin{cases}
t^{(m+2)(1-\frac{\sigma}{2\delta})+\frac{\sigma-2\beta}{2\delta}} \text{ if }n=2m+1, \\
t^{(m+1)(1-\frac{\sigma}{2\delta})+\frac{\sigma-2\beta}{2\delta}} \text{ if }n=2m.
\end{cases} \label{l4.1.10}
\end{equation}
Here we also note that $|\xi| \in (0,1]$, that is, $r \in (0,t^{\frac{1}{2\delta}}]$ and $rt^{-\frac{1}{2\delta}}\le 1$ which are useful in our proof. Summarizing, from (\ref{l4.1.1}) and (\ref{l4.1.7}) to (\ref{l4.1.10}) the statements of Lemma \ref{bd4.1.1} are proved.
\end{proof}

\begin{nx}
\fontshape{n}
\selectfont
Let us explain the result for the case $n=1$. We explained the proofs to Lemma \ref{bd4.1.1} for $n \ge 2$ only. However, in the case $n=1$ we only carry out partial integration with no need of the support of the vector field $Xf(r)$ as we did in (\ref{l4.1.3}). Then, following the steps of our considerations for odd spatial dimensions we may conclude that the statements of this lemma also hold for $n=1$.
\end{nx}
\noindent Following the proof of Lemma \ref{bd4.1.1} we may prove the following $L^1$ estimate, too.
\bbd \label{bd4.1.2}
The following estimate holds in $\R^n$ for any $n\ge 1$:
$$\Big\|F^{-1}\Big(e^{-c_1 |\xi|^{2\delta}t} |\xi|^{2\beta} \cos (c_2 |\xi|^\sigma t) \chi(|\xi|)\Big)(t,\cdot)\Big\|_{L^1}
\lesssim \begin{cases}
1 \text{ for } t\in (0,1], &\\
t^{(2+[\frac{n}{2}])(1-\frac{\sigma}{2\delta})-\frac{\beta}{\delta}} \text{ for } t\in[1,\ity), &
\end{cases}$$
with $\beta \ge 0$. Here $c_1$ is a positive and $c_2 \ne 0$ is a real constant.
\ebd

\noindent Finally, we consider oscillating integrals with more complicated oscillations in the integrand. We are going to prove the following result.
\bbd \label{bd4.1.3}
The following estimate holds in $\R^n$ for any $n\ge 1$:
$$\Big\|F^{-1}\Big(e^{-c_1 |\xi|^{2\delta}t}|\xi|^{2\beta} \frac{\sin (c_2 |\xi|^\sigma f(|\xi|) t)}{|\xi|^\sigma f(|\xi|)}\chi(|\xi|)\Big)(t,\cdot)\Big\|_{L^1} \lesssim \begin{cases}
t \text{ for } t\in (0,1], &\\
t^{(2+[\frac{n}{2}])(1-\frac{\sigma}{2\delta})+\frac{\sigma-2\beta}{2\delta}} \text{ for } t\in[1,\ity), &
\end{cases}$$
where $$f(|\xi|)=\sqrt{1-\frac{1}{4}|\xi|^{4\delta- 2\sigma}}$$
and $\beta \ge 0$. Here $c_1$ is a positive and $c_2 \ne 0$ is a real constant.
\ebd
\begin{proof}
We will follow the proof of Lemma $3.3$ in \cite{DaoReissig}. Hence, it is reasonable to present only the steps which are different. Then, we shall repeat some of the arguments as we did in the proof to Lemma \ref{bd4.1.1} to conclude the desired estimates.\medskip

\noindent First, let us consider $|x| \ge 1$ and $t \in (0,1]$. To obtain the first desired estimates in both cases of odd spatial dimensions $n=2m+1$ and even spatial dimensions $n=2m$ with $m \ge 1$, we assert the following estimates on the support of $\chi(r)$ and on the support of its derivatives:
$$ \Big| \partial_t^k \Big( \frac{\sin \big(c_2 r^\sigma f(r) t\big)}{f(r)} \chi(r) \Big) \Big| \lesssim r^{\sigma-k}t \text{ for all } k= 1,\cdots,m, $$
where
$$f(r)=\sqrt{1-\frac{1}{4}r^{4\delta- 2\sigma}}. $$
Here Fa\`{a} di Bruno's formula comes into play for all our estimates. We split the proof of the above estimate into several sub-steps as follows: \medskip

\noindent Step 1: $\quad$ Applying Proposition \ref{FadiBruno'sformula1} with $h(s)= \sqrt{s}$ and $g(r)=1-\frac{1}{4}r^{2(2\delta- \sigma)}$ we have
\begin{align*}
&\big|\partial_r^k f(r)\big| \lesssim \Big| \sum_{\substack{1\cdot m_1+\cdots+k\cdot m_k=k,\,m_i \ge 0}} g(r)^{\frac{1}{2}-(m_1+\cdots+m_k)}\prod_{j=1}^k \Big(-\frac{1}{4}r^{2(2\delta- \sigma)-j}\Big)^{m_j}\Big|\\
&\qquad \lesssim \sum_{\substack{1\cdot m_1+\cdots+k\cdot m_k=k,\, m_i \ge 0}} r^{2(2\delta- \sigma)(m_1+\cdots+m_k)-k}\lesssim r^{-k} \,\,\Big(\text{since}\,\,\frac{3}{4} \le g(r) \le 1 \,\,\,\text{for } \,r \le 1 \Big).
\end{align*}
In the same way we derive
\begin{equation}
\Big|\partial_r^k \Big(\frac{1}{f(r)} \Big)\Big| \lesssim r^{-k} \text{ for }k=1,\cdots,m. \label{l4.3.1}
\end{equation}
\noindent Step 2: $\quad$ Applying Proposition \ref{FadiBruno'sformula1} with $h(s)= \sin(c_2\,s)$ and $g(r)= r^\sigma f(r) t$ we get
\begin{align}
&\big| \partial_r^k \sin \big(c_2 r^\sigma f(r) t\big) \big| \lesssim \Big| \sum_{\substack{1\cdot m_1+\cdots+k\cdot m_k=k,\, m_i \ge 0}} \sin \big(c_2 r^\sigma f(r) t\big)^{(m_1+\cdots+m_k)} \prod_{j=1}^k \Big( \partial_r^j \big(r^\sigma f(r) t\big) \Big)^{m_j} \Big| \nonumber \\
&\qquad \lesssim \Big| \sum_{\substack{1\cdot m_1+\cdots+k\cdot m_k=k,\, m_i \ge 0}} \,\, \prod_{j=1}^k \Big( t\, \sum_{l=0}^j C_j^l r^{\sigma-j+l} f(r)^{(l)} \Big)^{m_j} \Big| \nonumber \\
&\qquad \lesssim \sum_{\substack{1\cdot m_1+\cdots+k\cdot m_k=k,\, m_i \ge 0}} \,\, \prod_{j=1}^k (t\, r^{\sigma-j})^{m_j} \lesssim \sum_{\substack{1\cdot m_1+\cdots+k\cdot m_k=k,\, m_i \ge 0}} r^{-k} (t\, r^\sigma)^{m_1+\cdots+m_k} \lesssim r^{\sigma-k} t. \label{l4.3.2}
\end{align}
Therefore, from (\ref{l4.3.1}) and (\ref{l4.3.2}) using the product rule for higher derivatives we may conclude
$$\Big| \partial_r^k \Big(\frac{\sin \big(c_2 r^\sigma f(r) t\big)}{f(r)}\Big) \Big| \lesssim  r^{\sigma-k} t \,\,\text{ for } \,\,k=1,\cdots,m. $$

\noindent Next, let us turn to consider $|x| \ge 1$ and $t \in [1,\ity)$. To derive the desired estimates by using similar ideas as in the proof to Lemma \ref{bd4.1.1}, we shall prove the following auxiliary estimates on the support of $\chi(t^{-\frac{1}{2\delta}}r)$ and on the support of its partial derivatives:
$$ \Big| \partial_r^k \Big( \frac{\sin \big(c_2 r^\sigma f(r) t^{1-\frac{\sigma}{2\delta}}\big)}{f(r)} \Big) \Big| \lesssim t^{1-\frac{\sigma}{2\delta}}r^{\sigma- k}(1+r^\sigma t^{1-\frac{\sigma}{2\delta}})^{k-1} \text{ if } k=1,\cdots,m, $$
where
$$f(r)= \sqrt{1-\frac{1}{4}t^{\frac{\sigma-2\delta}{\delta}}r^{2(2\delta- \sigma)}}. $$

\noindent Step 1: $\quad$ Applying Proposition \ref{FadiBruno'sformula1} with $h(s)=\sqrt{s}$ and $g(r)=1-\frac{1}{4}t^{\frac{\sigma-2\delta}{\delta}} r^{2(2\delta- \sigma)}$ we get
\begin{align*}
&\big|\partial_r^k f(r)\big| \lesssim \Big| \sum_{\substack{1\cdot m_1+\cdots+k\cdot m_k=k,\, m_i \ge 0}} g(r)^{\frac{1}{2}-(m_1+\cdots+m_k)}\prod_{j=1}^k \Big(-\frac{1}{4}t^{\frac{\sigma-2\delta}{\delta}}r^{2(2\delta- \sigma)-j}\Big)^{m_j}\Big|\\
&\qquad \lesssim \sum_{\substack{1\cdot m_1+\cdots+k\cdot m_k=k,\, m_i \ge 0}} \big( t^{\frac{\sigma-2\delta}{\delta}} r^{2(2\delta- \sigma)}\big)^{m_1+\cdots+m_k} r^{-k} \,\,\,  \Big(\text{since } \frac{3}{4} \le g(r) \le 1 \text{ for } r \le t^{\frac{1}{2\delta}}\Big)\\
&\qquad \lesssim r^{-k} \sum_{\substack{1\cdot m_1+\cdots+k\cdot m_k=k,\, m_i \ge 0}} (t^{-\frac{1}{2\delta}}r)^{2(2\delta- \sigma)(m_1+\cdots+m_k)} \lesssim r^{-k} \,\,\,  \Big(\text{since } t^{-\frac{1}{2\delta}}r \le 1 \text{ for } r \le t^{\frac{1}{2\delta}}\Big).
\end{align*}
An analogous treatment leads to
\begin{equation}
\Big| \partial_r^k \Big(\frac{1}{f(r)} \Big)\Big|  \lesssim r^{-k} \text{ for } k=1,\cdots,m. \label{l4.3.3}
\end{equation}

\noindent Step 2: $\quad$ Repeating the proof as we did in Lemma $3.3$ in \cite{DaoReissig} we have the following estimates:
\begin{equation}
\big| \partial_r^k \sin \big(c_2 r^\sigma f(r) t^{1-\frac{\sigma}{2\delta}}\big) \big| \lesssim t^{1-\frac{\sigma}{2\delta}} r^{\sigma-k} \big( 1+t^{1-\frac{\sigma}{2\delta}} r^\sigma\big)^{k-1}. \label{l4.3.4}
\end{equation}
Therefore, from (\ref{l4.3.3}) and (\ref{l4.3.4}) using the product rule for higher derivatives we may conclude
$$\Big| \partial_r^k \Big(\frac{\sin \big(c_2 r^\sigma f(r) t^{1-\frac{\sigma}{2\delta}}\big)}{f(r)} \Big)\Big| \lesssim t^{1-\frac{\sigma}{2\delta}} r^{\sigma-k} \big( 1+t^{1-\frac{\sigma}{2\delta}} r^\sigma\big)^{k-1} \text{ for } k=1,\cdots,m. $$
Summarizing, Lemma \ref{bd4.1.3} is proved.
\end{proof}

\noindent Following the steps of the proof to Lemma \ref{bd4.1.3} we may conclude the following statement, too.
\bbd \label{bd4.1.4}
The following estimate holds in $\R^n$ for any $n\ge 1$:
$$\Big\|F^{-1}\Big(e^{-c_1 |\xi|^{2\delta}t} |\xi|^{2\beta}\cos \big(c_2 |\xi|^\sigma f(|\xi|) t\big) \chi(|\xi|)\Big)(t,\cdot)\Big\|_{L^1}
\lesssim \begin{cases}
1 \text{ for } t\in (0,1], &\\
t^{(2+[\frac{n}{2}])(1-\frac{\sigma}{2\delta})- \frac{\beta}{\delta}} \text{ for } t\in[1,\ity), &
\end{cases}$$
where $$f(|\xi|)=\sqrt{1-\frac{1}{4}|\xi|^{4\delta- 2\sigma}}$$
and $\beta \ge 0$. Here $c_1$ is a positive and $c_2 \ne 0$ is a real constant.
\ebd

\begin{nx} \label{nx4.1.1}
\fontshape{n}
\selectfont
Following the proof of Lemmas from \ref{bd4.1.1} to \ref{bd4.1.4} we can conclude that all the desired statements still hold in the case $\delta= \sigma$.
\end{nx}

\renewcommand{\proofname}{Proof of Proposition \ref{md4.1.2}.}
\begin{proof}
In order to prove the first statement, by the relation (\ref{p3.2.1}) we replace $2\beta= a+2\delta$ and $2\beta= a$, respectively, in Lemmas \ref{bd4.1.3} and \ref{bd4.1.4}. For the sake of the relation (\ref{p3.2.2}), plugging $2\beta= a$ in Lemma \ref{bd4.1.3} we may conclude the second statement.
\end{proof}

\begin{nx} \label{nx4.1.2}
\fontshape{n}
\selectfont
Repeating the same arguments as in the proof of Proposition \ref{md4.1.2} we may see that these statements still hold in the case $\delta= \sigma$.
\end{nx}

\subsubsection{$L^1$ estimates for large frequencies} \label{Sec3.1.2}
\noindent Our approach is based on the paper \cite{NarazakiReissig}. According to the considerations in Section $5.2$ in \cite{NarazakiReissig}, with minor modifications in the steps of the proofs we may conclude the following $L^1$ estimates for large frequencies.

\bmd \label{md4.1.4}
The following estimates hold in $\mathbb{R}^n$ for any $n \geq 1$:
\begin{align*}
&\big\|F^{-1}\big(|\xi|^a \hat{K_0}\big(1-\chi(|\xi|)\big)\big)(t,\cdot)\big\|_{L^1}\lesssim
\begin{cases}
t^{-\frac{a}{2(\sigma-\delta)}} \text{ for } t\in (0,1], \\
t^{-\frac{a}{2\delta}} \text{ for } t\in[1,\ity), 
\end{cases} \\
&\big\|F^{-1}\big(|\xi|^a \hat{K_1}\big(1-\chi(|\xi|)\big)\big)(t,\cdot)\big\|_{L^1}\lesssim
\begin{cases}
t^{1-\frac{a}{2\delta}} \text{ for } t\in (0,1], \\
t^{1-\frac{a}{2(\sigma-\delta)}} \text{ for } t\in[1,\ity), 
\end{cases}
\end{align*}
for any non-negative number $a$.
\emd

\noindent Finally, from the statements of Propositions \ref{md4.1.2} and \ref{md4.1.4}, we may conclude the following $L^1$ estimates.
\bmd \label{md4.1.6}
The following estimates hold in $\R^n$ for any $n\ge 1$:
\begin{align*}
&\big\| F^{-1}\big(|\xi|^a \hat{K_0}\big)(t,\cdot)\big\|_{L^1}\lesssim \begin{cases}
t^{-\frac{a}{2(\sigma-\delta)}} \text{ for } t\in (0,1], &\\
t^{(2+[\frac{n}{2}])(1-\frac{\sigma}{2\delta})-\frac{a}{2\delta}} \text{ for } t\in[1,\ity), &
\end{cases}\\
&\big\| F^{-1}\big(|\xi|^a \hat{K_1}\big)(t,\cdot)\big\|_{L^1}\lesssim \begin{cases}
t^{1-\frac{a}{2\delta}} \text{ for } t\in (0,1], &\\
t^{1+(1+[\frac{n}{2}])(1-\frac{\sigma}{2\delta})-\frac{a}{2\delta}} \text{ for } t\in[1,\ity), &
\end{cases}
\end{align*}
for any non-negative number $a$.
\emd

\subsubsection{$L^\infty$ estimates} \label{Sec3.1.3}
\bmd \label{md4.1.9}
The following estimates hold in $\R^n$ for any $n\ge 1$:
\begin{align}
&\big\|F^{-1}\big(|\xi|^a \hat{K_0}\chi(|\xi|)\big)(t,\cdot)\big\|_{L^\ity}\lesssim
\begin{cases}
1 \text{ for } t \in (0,1], &\\
 t^{-\frac{n+a}{2\delta}} \text{ for } t \in [1,\ity), &
\end{cases} \label{p4.1.9.1} \\
&\big\|F^{-1}\big(|\xi|^a \hat{K_0}\big(1-\chi(|\xi|)\big)\big)(t,\cdot)\big\|_{L^\ity}\lesssim t^{-\frac{n+a}{2(\sigma-\delta)}} \text{ for } t \in (0,\ity), \label{p4.1.9.3} \\
&\big\|F^{-1}\big(|\xi|^a \hat{K_1}\chi(|\xi|)\big)(t,\cdot)\big\|_{L^\ity}\lesssim
\begin{cases}
t \text{ for } t \in (0,1], &\\
t^{1-\frac{n+a}{2\delta}} \text{ for } t \in [1,\ity), &
\end{cases} \label{p4.1.9.2} \\
&\big\|F^{-1}\big(|\xi|^a \hat{K_1}\big(1-\chi(|\xi|)\big)\big)(t,\cdot)\big\|_{L^\ity}\lesssim t^{1-\frac{n+a}{2(\sigma-\delta)}} \text{ for } t \in (0,\ity), \label{p4.1.9.4}
\end{align}
for any non-negative number $a$.
\emd

\renewcommand{\proofname}{Proof.}
\begin{proof}
Taking account of the representation for $\hat{K_1}$ we can re-write it as follows:
$$ \hat{K_1}(t,\xi)= e^{\lambda_1 t}\f{1- e^{(\lambda_2- \lambda_1)t}}{\lambda_1-\lambda_2}=
\begin{cases}
te^{\lambda_1 t}\int_0^1 e^{-\theta i \sqrt{4|\xi|^{2\sigma}-|\xi|^{4\delta}}t}d\theta \text{ for small } |\xi|, &\\
te^{\lambda_1 t}\int_0^1 e^{-\theta \sqrt{|\xi|^{4\delta}-4|\xi|^{2\sigma}}t}d\theta \text{ for large } |\xi|. &
\end{cases} $$
Thanks to the asymptotic behavior of the characteristic roots in (\ref{pt3.3}) and (\ref{pt3.4}), we arrive at
\begin{align*}
&\big|\hat{K_1}(t,\xi)\big| \lesssim t e^{-|\xi|^{2\delta}t},\,\, \big|\hat{K_0}(t,\xi) \big| \lesssim e^{-|\xi|^{2\delta}t} \text{ for small } |\xi|, \\
&\big|\hat{K_1}(t,\xi)\big| \lesssim t e^{-|\xi|^{2(\sigma- \delta)}t},\,\, \big|\hat{K_0}(t,\xi)\big| \lesssim e^{-c |\xi|^{2(\sigma- \delta)}t} \text{ for large } |\xi|,
\end{align*}
where $c$ is a suitable positive constant. Therefore, we may conclude all the statements that we wanted to prove.
\end{proof}

\begin{nx} \label{nx4.1.9}
\fontshape{n}
\selectfont
Following the approach to prove Proposition \ref{md4.1.9} we may notice that the statements (\ref{p4.1.9.1}) and (\ref{p4.1.9.2}) still hold in the case $\delta= \sigma$.
\end{nx}

\noindent From Proposition \ref{md4.1.9} the following statement follows immediately.
\bmd \label{md4.1.10}
The following estimates hold in $\R^n$ for any $n\ge 1$:
$$\big\|F^{-1}\big(|\xi|^a \hat{K_0}\big)(t,\cdot)\big\|_{L^\ity}\lesssim \begin{cases}
t^{-\frac{n+a}{2(\sigma-\delta)}} \text{ for } t\in (0,1], &\\
t^{-\frac{n+a}{2\delta}} \text{ for } t\in[1,\ity), &
\end{cases}$$
and
$$\big\|F^{-1}\big(|\xi|^a \hat{K_1}\big)(t,\cdot)\big\|_{L^\ity}\lesssim \begin{cases}
t^{1-\frac{n+a}{2(\sigma-\delta)}} \text{ for } t\in (0,1], &\\
t^{1-\frac{n+a}{2\delta}} \text{ for } t\in[1,\ity), &
\end{cases}$$
for any non-negative number $a$.
\emd

\subsubsection{$L^r$ estimates}\label{Sec3.1.4}

\noindent From the statements of Propositions \ref{md4.1.6} and \ref{md4.1.10}, by employing an interpolation argument we may conclude the following $L^r$ estimates.
\bmd \label{md4.1.12}
The following estimates hold in $\R^n$ for $n\ge 1$:
\begin{align*}
\big\|F^{-1}\big(|\xi|^a \hat{K_0}\big)(t,\cdot)\big\|_{L^r}&\lesssim
\begin{cases}
t^{-\frac{n}{2(\sigma-\delta)}(1-\frac{1}{r})-\frac{a}{2(\sigma-\delta)}} \text{ if } t\in (0,1], &\\
t^{(2+[\frac{n}{2}])(1-\frac{\sigma}{2\delta})\frac{1}{r} -\frac{n}{2\delta}(1-\frac{1}{r})-\frac{a}{2\delta}} \text{ if } t\in[1,\ity), &
\end{cases} \\
\big\|F^{-1}\big(|\xi|^a \hat{K_1}\big)(t,\cdot)\big\|_{L^r}&\lesssim
\begin{cases}
t^{1-\frac{n}{2(\sigma-\delta)}(1-\frac{1}{r})-\frac{a}{2(\sigma-\delta)}} \text{ if } t\in (0,1], &\\
t^{1+(1+[\frac{n}{2}])(1-\frac{\sigma}{2\delta})\frac{1}{r} -\frac{n}{2\delta}(1-\frac{1}{r})-\frac{a}{2\delta}} \text{ if } t\in[1,\ity), &
\end{cases}
\end{align*}
for all $r \in [1,\ity]$ and any non-negative number $a$.
\emd

\bhq{$L^p- L^q$ estimates not necessarily on the conjugate line} \label{hq4.3.1}\\
Let $\delta\in \big(\frac{\sigma}{2},\sigma\big)$ in (\ref{pt6.3}) and $1\le p\le q\le \ity$. Then, the solutions to (\ref{pt6.3}) satisfy the $L^p- L^q$ estimates
\begin{align*}
\big\||D|^a u(t,\cdot )\big\|_{L^q}& \lesssim
\begin{cases}
t^{-\frac{n}{2(\sigma-\delta)}(1-\frac{1}{r})-\frac{a}{2(\sigma-\delta)}} \|u_0\|_{L^p}+ t^{1-\frac{n}{2(\sigma-\delta)}(1-\frac{1}{r})-\frac{a}{2(\sigma-\delta)}} \|u_1\|_{L^p} \text{ if } t\in (0,1], &\\
t^{(2+[\frac{n}{2}])(1-\frac{\sigma}{2\delta})\frac{1}{r} -\frac{n}{2\delta}(1-\frac{1}{r})-\frac{a}{2\delta}} \|u_0\|_{L^p}\\
\qquad+ t^{1+(1+[\frac{n}{2}])(1-\frac{\sigma}{2\delta})\frac{1}{r} -\frac{n}{2\delta}(1-\frac{1}{r})-\frac{a}{2\delta}} \|u_1\|_{L^p} \text{ if } t\in[1,\ity), &
\end{cases}\\
\big\||D|^a u_t(t,\cdot )\big\|_{L^q}& \lesssim
\begin{cases}
t^{-\frac{n}{2(\sigma-\delta)}(1-\frac{1}{r})-\frac{a+2\delta}{2(\sigma-\delta)}} \|u_0\|_{L^p}+ t^{1-\frac{n}{2(\sigma-\delta)}(1-\frac{1}{r})-\frac{a+2\delta}{2(\sigma-\delta)}} \|u_1\|_{L^p} \text{ if } t\in (0,1], &\\
t^{1+(1+[\frac{n}{2}])(1-\frac{\sigma}{2\delta})\frac{1}{r} -\frac{n}{2\delta}(1-\frac{1}{r})-\frac{a+2{\sigma-\delta}}{2\delta}} \|u_0\|_{L^p}\\
\qquad+ t^{(2+[\frac{n}{2}])(1-\frac{\sigma}{2\delta})\frac{1}{r} -\frac{n}{2\delta}(1-\frac{1}{r})- \frac{a}{2\delta}} \|u_1\|_{L^p} \text{ if } t\in[1,\ity), &
\end{cases}
\end{align*}
where $1+ \frac{1}{q}= \frac{1}{r}+\frac{1}{p}$, for any non-negative number $a$ and for all $n\ge 1$.
\ehq

\begin{proof}
For the sake of the statements in Proposition \ref{md4.1.12}, applying Young's convolution inequality we may conclude the first statement. In order to prove some estimates for the time derivative of the solution we notice that the following relations hold:
$$ \partial_t \hat{K_0}= -|\xi|^{2\sigma}\hat{K_1}\text{ and } \partial_t \hat{K_1}= \hat{K_0}- |\xi|^{2\delta}\hat{K_1}. $$
Hence, employing again Young's convolution inequality and using Proposition \ref{md4.1.12} we may conclude the second statement. Summarizing, Corollary \ref{hq4.3.1} is proved.
\end{proof}

\begin{nx} \label{Rem3.5}
\fontshape{n}
\selectfont
Let us compare our results with some known results in \cite{NarazakiReissig}. In the special case $\sigma=1$ one may observe that the decay estimates for the solution itself appearing in Corollary \ref{hq4.3.1} are asymptotically the same as the corresponding ones obtained in the cited paper if we consider the case of sufficiently large space dimensions $n$.
\end{nx}

\begin{nx} \label{Rem3.6}
\fontshape{n}
\selectfont
We can see that there appears the singular behavior of the time-dependent coefficients for $t \to +0$ in Corollary \ref{hq4.3.1}. This causes some difficulties to treat the semi-linear models. Hence, in order to overcome this, we state the following corollary by assuming additional regularity for the data.
\end{nx}

\bhq {$L^{q}- L^{q}$ estimates with additional $L^{m}$ regularity for the data} \label{hq4.3.2} \\
Let $\delta\in \big(\f{\sigma}{2},\sigma \big)$ in (\ref{pt6.3}), $q\in (1,\ity)$ and $m\in [1,q)$. Then the Sobolev solutions to (\ref{pt6.3}) satisfy the following $(L^m \cap L^q)-L^q$ estimates:
\begin{align*}
\big\||D|^a u(t,\cdot)\big\|_{L^q}& \lesssim  (1+t)^{(2+[\frac{n}{2}])(1-\frac{\sigma}{2\delta})\frac{1}{r} -\frac{n}{2\delta}(1-\frac{1}{r})-\frac{a}{2\delta}} \|u_0\|_{L^m \cap H^{a,q}} \nonumber\\
& \qquad \quad + (1+t)^{1+(1+[\frac{n}{2}])(1-\frac{\sigma}{2\delta})\frac{1}{r} -\frac{n}{2\delta}(1-\frac{1}{r})-\frac{a}{2\delta}}\|u_1\|_{L^m \cap H^{[a-2\delta]^+,q}}, \\
\big\||D|^a u_t(t,\cdot)\big\|_{L^q}& \lesssim  (1+t)^{(1+[\frac{n}{2}])(1-\frac{\sigma}{2\delta})\frac{1}{r} -\frac{n}{2\delta}(1-\frac{1}{r})-\frac{a+2(\sigma-\delta)}{2\delta}} \|u_0\|_{L^m \cap H^{a+2(\sigma-\delta),q}} \nonumber\\\
& \qquad \quad + (1+t)^{(2+[\frac{n}{2}])(1-\frac{\sigma}{2\delta})\frac{1}{r} -\frac{n}{2\delta}(1-\frac{1}{r})-\frac{a}{2\delta}}\|u_1\|_{L^m \cap H^{a,q}}.
\end{align*}
Moreover, the following $L^q-L^q$ estimates are satisfied:
\begin{align*}
\big\||D|^a u(t,\cdot)\big\|_{L^q}& \lesssim (1+t)^{(2+[\frac{n}{2}])(1-\frac{\sigma}{2\delta})-\frac{a}{2\delta}} \|u_0\|_{H^{a,q}}+ (1+t)^{1+(1+[\frac{n}{2}])(1-\frac{\sigma}{2\delta})-\frac{a}{2\delta}} \|u_1\|_{H^{[a-2\delta]^+,q}}, \\
\big\||D|^a u_t(t,\cdot)\big\|_{L^q}& \lesssim (1+t)^{(1+[\frac{n}{2}])(1-\frac{\sigma}{2\delta})-\frac{a+2(\sigma-\delta)}{2\delta}} \|u_0\|_{H^{a+2(\sigma-\delta),q}}+ (1+t)^{(2+[\frac{n}{2}])(1-\frac{\sigma}{2\delta})-\frac{a}{2\delta}} \|u_1\|_{H^{a,q}}.
\end{align*}
Here $1+ \frac{1}{q}= \frac{1}{r}+\frac{1}{m}$, $a$ is a non-negative number and the dimension $n\ge 1$.
\ehq

\begin{proof}
To derive the $(L^m \cap L^q)- L^q$ estimates, on the one hand we control the $L^q$ norm of the small frequency part of the solutions by the $L^m$ norm of the data. On the other hand, its high-frequency part is estimated by using the $L^q-L^q$ estimates with a suitable regularity of the data depending on the order $a$ of derivatives. Finally, applying Young's convolution inequality we may conclude all the statements what we wanted to prove.
\end{proof}

\subsection{Estimates of oscillating integrals in the case of visco-elastic type damping: $\delta=\sigma$} \label{Sec3.2}

First, let us explain our strategies to deal with estimates in the case $\delta= \sigma$. We recall that in the case $\delta \in (\frac{\sigma}{2},\sigma)$ our goal is to obtain $L^p-L^q$ estimates not necessarily on the conjugate line with $1\le p\le q\le \ity$. For this reason, we need to develop some techniques from the paper \cite{NarazakiReissig} in order to conclude $L^1$ estimates, $L^\ity$ estimates and $L^r$ estimates, with $r \in [1,\ity]$, as well. Moreover, this strategy is also applied with an extension to the case $\delta= \sigma$ to get these estimates for small frequencies (see later, Section \ref{Sec3.2.1}). Meanwhile, for large frequencies in the case $\delta= \sigma$ this strategy fails to derive $L^1$ estimates as in Proposition \ref{md4.1.4}, and $L^\ity$ estimates as in Proposition \ref{md4.1.9} for (\ref{p4.1.9.3}) and (\ref{p4.1.9.4}). Hence, it is reasonable to apply the Mikhlin- H\"{o}mander multiplier theorem for large frequencies in the case $\delta= \sigma$. It is clear that this approach is only to conclude $L^q-L^q$ estimates for large frequencies with $q \in (1,\ity)$. Therefore, in the case $\delta= \sigma$ the aim to obtain $L^p-L^q$ estimates not necessarily on the conjugate line with $1\le p\le q\le \ity$ is beyond the scope of our paper.

\subsubsection{$L^p- L^q$ estimates not necessarily on the conjugate line for small frequencies}\label{Sec3.2.1}

\noindent By interpolation argument, from Remarks \ref{nx4.1.2} and \ref{nx4.1.9} we can conclude the following $L^r$ estimates for small frequencies.
\bmd \label{md0.1.3}
The following estimates hold in $\R^n$ for any $n\ge 1$:
\begin{align*}
\big\|F^{-1}\big(|\xi|^a \hat{K_0}\chi(|\xi|)\big)(t,\cdot)\big\|_{L^r}&\lesssim
\begin{cases}
1 \text{ if } t\in (0,1], &\\
t^{\frac{1}{2}(2+[\frac{n}{2}])\frac{1}{r} -\frac{n}{2\sigma}(1-\frac{1}{r})-\frac{a}{2\sigma}} \text{ if } t\in[1,\ity), &
\end{cases} \\
\big\|F^{-1}\big(|\xi|^a \hat{K_1}\chi(|\xi|)\big)(t,\cdot)\big\|_{L^r}&\lesssim
\begin{cases}
t \text{ if } t\in (0,1], &\\
t^{1+\frac{1}{2}(1+[\frac{n}{2}])\frac{1}{r} -\frac{n}{2\sigma}(1-\frac{1}{r})-\frac{a}{2\sigma}} \text{ if } t\in[1,\ity), &
\end{cases}
\end{align*}
for all $r \in [1,\ity]$ and any non-negative number $a$.
\emd

Repeating the proof of Corollary \ref{hq4.3.1} we obtain the following estimate by using Proposition \ref{md0.1.3}.
\bhq \label{hq4.3.3}
Let $\delta=\sigma$ in (\ref{pt6.3}) and $1\le p\le q\le \ity$. Then, the solutions to (\ref{pt6.3}) satisfy the $L^p- L^q$ estimates
\begin{align*}
\big\||D|^a u_{\chi}(t,\cdot )\big\|_{L^q}& \lesssim (1+t)^{\frac{1}{2}(2+[\frac{n}{2}])\frac{1}{r} -\frac{n}{2\sigma}(1-\frac{1}{r})-\frac{a}{2\sigma}} \|u_0\|_{L^p}+ t(1+t)^{\frac{1}{2}(1+[\frac{n}{2}])\frac{1}{r} -\frac{n}{2\sigma}(1-\frac{1}{r})-\frac{a}{2\sigma}} \|u_1\|_{L^p}, \\
\big\||D|^a \partial_t u_{\chi}(t,\cdot )\big\|_{L^q}& \lesssim (1+t)^{\frac{1}{2}(1+[\frac{n}{2}])\frac{1}{r} -\frac{n}{2\sigma}(1-\frac{1}{r})-\frac{a}{2\sigma}} \|u_0\|_{L^p}+ (1+t)^{\frac{1}{2}(2+[\frac{n}{2}])\frac{1}{r} -\frac{n}{2\sigma}(1-\frac{1}{r})-\frac{a}{2\sigma}} \|u_1\|_{L^p},
\end{align*}
where $1+ \frac{1}{q}= \frac{1}{r}+\frac{1}{p}$, for any non-negative number $a$ and for all $n\ge 1$.
\ehq

\begin{nx}
\fontshape{n}
\selectfont
Let us compare our results with some known results in \cite{Shibata}. In the special case  $\sigma=\delta=1$ we may observe that the time-dependent coefficients in the $L^1-L^1$ estimate for the solution itself appearing in Corollary \ref{hq4.3.3} are asymptotically the same as the corresponding ones obtained in the cited paper if we consider the case of sufficiently large space dimensions $n$.
\end{nx}

\subsubsection{$L^q- L^q$ estimates for large frequencies} \label{Sec3.2.2}
\noindent First, we can re-write the characteristic roots as follows:
\begin{equation}
\lambda_1(\xi)= -1-\mu(\xi) \text{ and } \lambda_2(\xi)= -|\xi|^{2\sigma}+1+\mu(\xi), \label{characteristicRoot1}
\end{equation}
where
\begin{equation}
\mu(\xi)=-1+g\Big(\frac{4}{|\xi|^{2\sigma}}\Big) \text{ and } g(s)=\int_0^{1}(1-\theta s)^{-\frac{1}{2}}d\theta. \label{characteristicRoot2}
\end{equation}
Now, we introduce the following abbreviations:
\begin{align*}
&K_0^1=K_0^{1}(t,x):= F^{-1}\Big(\frac{\lambda_2(\xi) e^{\lambda_1(\xi) t}}{\lambda_1(\xi)- \lambda_2(\xi)}v_0(\xi)\big(1-\chi(\xi)\big) \Big)(t,x),\\
&K_0^2=K_0^{2}(t,x):= F^{-1}\Big(\frac{\lambda_1(\xi) e^{\lambda_2(\xi) t}}{\lambda_1(\xi)- \lambda_2(\xi)}v_0(\xi)\big(1-\chi(\xi)\big) \Big)(t,x),\\
&K_1^1=K_1^{1}(t,x):= F^{-1}\Big(\frac{e^{\lambda_1(\xi) t}}{\lambda_1(\xi)- \lambda_2(\xi)}v_1(\xi)\big(1-\chi(\xi)\big) \Big)(t,x),\\
&K_1^2=K_1^{2}(t,x):= F^{-1}\Big(\frac{e^{\lambda_2(\xi) t}}{\lambda_1(\xi)- \lambda_2(\xi)}v_1(\xi)\big(1-\chi(\xi)\big) \Big)(t,x).
\end{align*}
We shall prove the following results.
\bmd \label{md0.2.1}
Let $q\in (1,\ity)$. Then, the following estimates hold:
\begin{align*}
&\big\|\partial_t^{j}|D|^a K_0^{1}(t,\cdot)\big\|_{L^q} \lesssim e^{-ct} \|u_0\|_{H^{a,q}}, \\
&\big\|\partial_t^{j}|D|^a K_0^{2}(t,\cdot)\big\|_{L^q} \lesssim e^{-ct} \|u_0\|_{H^{[2\sigma j-2\sigma+a]^{+},q}}, \\
&\big\|\partial_t^{j}|D|^a K_1^{1}(t,\cdot)\big\|_{L^q} \lesssim e^{-ct} \|u_1\|_{H^{[a-2\sigma]^{+},q}}, \\
&\big\|\partial_t^{j}|D|^a K_1^{2}(t,\cdot)\big\|_{L^q} \lesssim e^{-ct} \|u_1\|_{H^{[2\sigma j-2\sigma+a]^{+},q}},
\end{align*}
for any $t> 0$, $a \ge 0$, integer $j \ge 0$ and a suitable positive constant $c$.
\emd

According to application of the Mikhlin- H\"{o}mander multiplier theorem (see also \cite{DabbiccoEbert,Miyachi}) for Fourier multipliers from Proposition \ref{PropositionMultiplier}, in order to prove Proposition \ref{md0.2.1} we shall show the following auxiliary estimates.
\bbd \label{bd0.2.1}
The following estimates hold in $\R^n$ for sufficiently large $|\xi|$:
\begin{align}
&\big|\partial_\xi^\alpha |\xi|^{-2\sigma}\big| \lesssim |\xi|^{-2\sigma-|\alpha|} \text{ for all }\alpha, \label{bd0.2.1.1} \\
&\big|\partial_\xi^\alpha |\xi|^{2p\sigma} \big| \lesssim |\xi|^{2p\sigma-|\alpha|}\text{ for all }\alpha\text{ and } p\in\R, \label{bd0.2.1.11} \\
&\Big|\partial_\xi^\alpha g\Big(\frac{4}{|\xi|^{2\sigma}}\Big)\Big| \lesssim |\xi|^{-2\sigma-|\alpha|}\text{ for all }|\alpha|\ge 1, \text{ and }\Big| g\Big(\frac{4}{|\xi|^{2\sigma}}\Big)\Big| \lesssim 1, \label{bd0.2.1.2} \\
&\big|\partial_\xi^\alpha \mu(\xi)\big| \lesssim |\xi|^{-2\sigma-|\alpha|}\text{ for all }\alpha, \label{bd0.2.1.3} \\
&\big|\partial_\xi^\alpha \lambda_2(\xi)\big| \lesssim |\xi|^{2\sigma-|\alpha|}\text{ for all }\alpha, \label{bd0.2.1.4} \\
&\big|\partial_\xi^\alpha \lambda_1(\xi)\big| \lesssim |\xi|^{-2\sigma-|\alpha|}\text{ for all }|\alpha|\ge 1, \text{ and }\big|\lambda_1(\xi)\big| \lesssim 1, \label{bd0.2.1.5} \\
&\Big|\partial_\xi^\alpha \big(\lambda_1(\xi)-\lambda_2(\xi)\big)^{-1}\Big| \lesssim |\xi|^{-2\sigma-|\alpha|}\text{ for all }\alpha, \label{bd0.2.1.6} \\
&\big|\partial_\xi^\alpha \lambda_2^j(\xi) \big| \lesssim |\xi|^{2\sigma j-|\alpha|}\text{ for all }\alpha \text{ and } j\ge 0, \label{bd0.2.1.7} \\
&\big|\partial_\xi^\alpha \lambda_1^j(\xi) \big| \lesssim |\xi|^{-|\alpha|}\text{ for all }\alpha \text{ and } j\ge 0, \label{bd0.2.1.8} \\
&\big|\partial_\xi^\alpha \big(|\xi|^b \lambda_2^j(\xi)\big) \big| \lesssim |\xi|^{2\sigma j+b-|\alpha|}\text{ for all }\alpha, \text{ for any } b\in \R \text{ and } j\ge 0, \label{bd0.2.1.9} \\
&\big|\partial_\xi^\alpha \big(|\xi|^b \lambda_1^j(\xi)\big) \big| \lesssim |\xi|^{b-|\alpha|}\text{ for all }\alpha, \text{ for any } b\in \R \text{ and } j\ge 0, \label{bd0.2.1.10} \\
&\big|\partial_\xi^\alpha \big(e^{\lambda_2(\xi)t}\big) \big| \lesssim e^{-ct}|\xi|^{-|\alpha|}, \label{bd0.2.1.12}\\
&\text{ for all }\alpha \text{ and } t>0, \text{ where $c$ is a suitable positive constant}, \nonumber \\
&\big|\partial_\xi^\alpha \big(e^{\lambda_1(\xi)t}\big) \big| \lesssim e^{-ct}|\xi|^{-|\alpha|}, \label{bd0.2.1.13} \\
&\text{ for all }\alpha \text{ and }t>0, \text{ where $c$ is a suitable positive constant}, \nonumber \\
&\Big|\partial_\xi^\alpha \Big(\frac{\lambda_1(\xi) e^{\lambda_2(\xi)t} \lambda_2^j(\xi) |\xi|^b}{\lambda_1(\xi)- \lambda_2(\xi)} \Big) \Big| \lesssim e^{-ct}|\xi|^{2\sigma j+b-2\sigma-|\alpha|}, \label{bd0.2.1.14} \\
&\text{ for all }\alpha, \text{ for any } b\in \R,\, j\ge 0 \text{ and } t>0, \text{ where $c$ is a suitable positive constant}, \nonumber \\
&\Big|\partial_\xi^\alpha \Big(\frac{e^{\lambda_2(\xi)t} \lambda_2^j(\xi) |\xi|^b}{\lambda_1(\xi)- \lambda_2(\xi)} \Big) \Big| \lesssim e^{-ct}|\xi|^{2\sigma j+b-2\sigma-|\alpha|}, \label{bd0.2.1.15} \\
&\text{ for all }\alpha, \text{ for any } b\in \R,\, j\ge 0 \text{ and } t>0, \text{ where $c$ is a suitable positive constant}, \nonumber \\
&\Big|\partial_\xi^\alpha \Big(\frac{\lambda_2(\xi) e^{\lambda_1(\xi)t} \lambda_1^j(\xi) |\xi|^b }{\lambda_1(\xi)- \lambda_2(\xi)} \Big) \Big| \lesssim e^{-ct}|\xi|^{b-|\alpha|}, \label{bd0.2.1.16} \\
&\text{ for all }\alpha, \text{ for any } b\in \R,\, j\ge 0 \text{ and } t>0, \text{ where $c$ is a suitable positive constant}, \nonumber \\
&\Big|\partial_\xi^\alpha\Big(\frac{e^{\lambda_1(\xi)t} \lambda_1^j(\xi) |\xi|^b}{\lambda_1(\xi)- \lambda_2(\xi)} \Big) \Big| \lesssim e^{-ct}|\xi|^{b-2\sigma-|\alpha|}, \label{bd0.2.1.17} \\
&\text{ for all }\alpha, \text{ for any } b\in \R,\, j\ge 0 \text{ and } t>0, \text{ where $c$ is a suitable positive constant}, \nonumber
\end{align}
\ebd

\begin{proof}
 In order to prove all statements in Lemma \ref{bd0.2.1}, we shall apply Lemma \ref{LemmaDerivative} and Leibniz rule of the multivariable calculus. Indeed, we will indicate the proof of the above estimates as follows:\medskip

\noindent First, we can that the proof of (\ref{bd0.2.1.1}) and (\ref{bd0.2.1.11}) is trivial. By (\ref{bd0.2.1.1}) applying Lemma \ref{LemmaDerivative} with $h(s)=g(s)$ and $f(\xi)=4|\xi|^{-2\sigma}$ we can conclude (\ref{bd0.2.1.2}). By (\ref{characteristicRoot1}) and (\ref{characteristicRoot2}), the statements from (\ref{bd0.2.1.3}) to (\ref{bd0.2.1.5}) are immediately follow from (\ref{bd0.2.1.11}) to (\ref{bd0.2.1.2}). In the analogous way, by (\ref{bd0.2.1.4}) and (\ref{bd0.2.1.5}) we get (\ref{bd0.2.1.6}) with $h(s)= s^{-1}$ and $f(\xi)=\lambda_1(\xi)-\lambda_2(\xi)$. By (\ref{bd0.2.1.4}) we obtain (\ref{bd0.2.1.7}) with $h(s)=s^j$ and $f(\xi)=\lambda_2(\xi)$. By (\ref{bd0.2.1.5}) we derive (\ref{bd0.2.1.8}) with $h(s)=s^j$ and $f(\xi)=\lambda_1(\xi)$. Using the Leibniz rule we conclude (\ref{bd0.2.1.9}) after using (\ref{bd0.2.1.11}) and (\ref{bd0.2.1.7}). Analogously, we obtain (\ref{bd0.2.1.10}) by using (\ref{bd0.2.1.11}) and (\ref{bd0.2.1.8}). Applying Lemma \ref{LemmaDerivative} with $h(s)=e^{st}$ and $f(\xi)=\lambda_2(\xi)$ we have (\ref{bd0.2.1.12}) by taking account of (\ref{bd0.2.1.4}), where we note that $\lambda_2(\xi) \le -\frac{1}{2}|\xi|^{2\sigma}$. In the same way, by (\ref{bd0.2.1.5}) and $\lambda_1(\xi) \le -1$ we get (\ref{bd0.2.1.13}) with $h(s)=e^{st}$ and $f(\xi)=\lambda_1(\xi)$. Combining (\ref{bd0.2.1.5}), (\ref{bd0.2.1.6}), (\ref{bd0.2.1.9}) and (\ref{bd0.2.1.12}) we may conclude (\ref{bd0.2.1.14}) and (\ref{bd0.2.1.15}) by using the Leibniz rule. Finally, combining (\ref{bd0.2.1.4}), (\ref{bd0.2.1.6}), (\ref{bd0.2.1.10}) and (\ref{bd0.2.1.13}) we arrive at (\ref{bd0.2.1.16}) and (\ref{bd0.2.1.17}) by using again the Leibniz rule.
\end{proof}

\renewcommand{\proofname}{Proof of Proposition \ref{md0.2.1}.}
\begin{proof}
First, taking account of estimates for $K_0^{2}$ and some its derivatives we will divide our consideration into two cases. In the first case, if $2\sigma j-2\sigma+a \ge 0$, then we can write
$$\partial_t^{j}|D|^a K_0^{2}(t,x)= F^{-1}\Big(\frac{ \lambda_1(\xi) e^{\lambda_2(\xi)t} \lambda_2^j(\xi) |\xi|^{2\sigma -2\sigma j}}{\lambda_1(\xi)- \lambda_2(\xi)}\big(1-\chi(\xi)\big) |\xi|^{2\sigma j-2\sigma+a}v_0(\xi) \Big)(t,x). $$
By choosing $b=2\sigma -2\sigma j$ in (\ref{bd0.2.1.14}), we get for all $\alpha$ the estimates
$$\Big|\partial_\xi^\alpha \Big(\frac{ \lambda_1(\xi) e^{\lambda_2(\xi)t} \lambda_2^j(\xi) |\xi|^{2\sigma -2\sigma j}}{\lambda_1(\xi)- \lambda_2(\xi)} \Big)\Big| \lesssim e^{-ct}|\xi|^{-|\alpha|}, $$
where $c$ is a suitable positive constant. By Proposition \ref{PropositionMultiplier}, we can conclude
\begin{equation}
\big\|\partial_t^{j}|D|^a K_0^{2}(t,\cdot)\big\|_{L^q} \lesssim e^{-ct} \|u_0\|_{H^{2\sigma j-2\sigma+a  ,q}}. \label{p0.2.1.1}
\end{equation}
In the second case, if $2\sigma j-2\sigma+a < 0$, then we can write
$$\partial_t^{j}|D|^a K_0^{2}(t,x)= F^{-1}\Big(\frac{ \lambda_1(\xi) e^{\lambda_2(\xi)t} \lambda_2^j(\xi) |\xi|^a}{\lambda_1(\xi)- \lambda_2(\xi)}\big(1-\chi(\xi)\big) v_0(\xi) \Big)(t,x). $$
By choosing $b=a$ in (\ref{bd0.2.1.14}), we derive for all $\alpha$ the estimates
$$\Big|\partial_\xi^\alpha \Big(\frac{\lambda_1(\xi) e^{\lambda_2(\xi)t} \lambda_2^j(\xi) |\xi|^a}{\lambda_1(\xi)- \lambda_2(\xi)}\Big)\Big| \lesssim e^{-ct}|\xi|^{2\sigma j+a-2\sigma-|\alpha|}\lesssim e^{-ct}|\xi|^{-|\alpha|}, $$
where $c$ is a suitable positive constant. By Proposition \ref{PropositionMultiplier}, we may conclude
\begin{equation}
\big\|\partial_t^{j}|D|^a K_0^{2}(t,\cdot)\big\|_{L^q} \lesssim e^{-ct} \|u_0\|_{L^q}. \label{p0.2.1.2}
\end{equation}
Hence, from (\ref{p0.2.1.1}) and (\ref{p0.2.1.2}) we have proved the second statement in Proposition \ref{md0.2.1}. By the same arguments we may also conclude the last statement in Proposition \ref{md0.2.1} by using (\ref{bd0.2.1.15}). Analogously, in order to estimate the third statement we will apply $b=2\sigma$ and $b=a$ in (\ref{bd0.2.1.17}), respectively, if $a \ge 2\sigma$ and $a < 2\sigma$. Finally, using (\ref{bd0.2.1.16}) with $b=0$ immediately leads to the remaining statement. Summarizing, the proof to Proposition \ref{md0.2.1} is completed.
\end{proof}

\begin{nx} \label{nx0.2.1}
\fontshape{n}
\selectfont
The exponential decay $e^{-ct}$  appearing in Proposition \ref{md0.2.1} for large frequencies is better than the potential decay in Proposition \ref{md0.1.3} for small frequencies. Since we have in mind that the real part of the characteristic roots $\lambda_{1,2}$ is negative in the middle zone $|\xi| \in \{\e,\, \frac{1}{\e}\}$ with sufficiently small $\e$, the corresponding estimates yield  in this zone an exponential decay, too.
\end{nx}

\noindent From Proposition \ref{md0.2.1} we conclude the following estimates for large frequencies.
\bhq \label{hq4.3.4}
Let $\delta=\sigma$ in (\ref{pt6.3}) and $q\in (1,\ity)$. Then, the solutions to (\ref{pt6.3}) satisfy the $L^q- L^q$ estimates
$$ \big\|\partial_t^{j}|D|^a u_{1-\chi}(t,\cdot)\big\|_{L^q} \lesssim e^{-ct}\big(\|(u_0,u_1)\|_{H^{[2\sigma j-2\sigma+a]^{+},q}}+\|u_0\|_{H^{a,q}}+\|u_1\|_{H^{[a-2\sigma]^{+},q}}\big), $$
for any $t> 0$, $a \ge 0$, integer $j\ge 0$ and a suitable positive constant $c$.
\ehq

\begin{nx}
\fontshape{n}
\selectfont
Let us compare our results with some known results in \cite{Shibata}. In the special case  $\sigma=\delta=1$ we may observe that the decay rates for $L^q-L^q$ estimate for the solution itself appearing in Corollary \ref{hq4.3.4} are exactly the same as the corresponding ones obtained in the cited paper.
\end{nx}

\noindent From Corollaries \ref{hq4.3.3} and \ref{hq4.3.4} we conclude the following estimates.
\bhq {$L^{q}- L^{q}$ estimates with additional $L^{m}$ regularity for the data} \label{hq4.3.5} \\
Let $\delta= \sigma$ in (\ref{pt6.3}), $q\in (1,\ity)$ and $m\in [1,q)$. Then the Sobolev solutions to (\ref{pt6.3}) satisfy the following $(L^m \cap L^q)-L^q$ estimates:
\begin{align*}
\big\||D|^a u(t,\cdot)\big\|_{L^q}& \lesssim  (1+t)^{\frac{1}{2}(2+[\frac{n}{2}])\frac{1}{r} -\frac{n}{2\sigma}(1-\frac{1}{r})-\frac{a}{2\sigma}} \|u_0\|_{L^m \cap H^{a,q}} \nonumber\\
& \qquad \quad + (1+t)^{1+\frac{1}{2}(1+[\frac{n}{2}])\frac{1}{r} -\frac{n}{2\sigma}(1-\frac{1}{r})-\frac{a}{2\sigma}}\|u_1\|_{L^m \cap H^{[a-2\sigma]^+,q}}, \\
\big\||D|^a u_t(t,\cdot)\big\|_{L^q}& \lesssim  (1+t)^{\frac{1}{2}(1+[\frac{n}{2}])\frac{1}{r} -\frac{n}{2\sigma}(1-\frac{1}{r})-\frac{a}{2\sigma}} \|u_0\|_{L^m \cap H^{a,q}} \nonumber\\\
& \qquad \quad + (1+t)^{\frac{1}{2}(2+[\frac{n}{2}])\frac{1}{r} -\frac{n}{2\sigma}(1-\frac{1}{r})-\frac{a}{2\sigma}}\|u_1\|_{L^m \cap H^{a,q}}.
\end{align*}
Moreover, the following $L^q-L^q$ estimates are satisfied:
\begin{align*}
\big\||D|^a u(t,\cdot)\big\|_{L^q}& \lesssim (1+t)^{\frac{1}{2}(2+[\frac{n}{2}])-\frac{a}{2\sigma}} \|u_0\|_{H^{a,q}}+ (1+t)^{1+\frac{1}{2}(1+[\frac{n}{2}])-\frac{a}{2\sigma}} \|u_1\|_{H^{[a-2\sigma]^+,q}}, \\
\big\||D|^a u_t(t,\cdot)\big\|_{L^q}& \lesssim (1+t)^{\frac{1}{2}(1+[\frac{n}{2}])-\frac{a}{2\sigma}} \|u_0\|_{H^{a,q}}+ (1+t)^{\frac{1}{2}(2+[\frac{n}{2}])-\frac{a}{2\sigma}} \|u_1\|_{H^{a,q}}.
\end{align*}
Here $1+ \frac{1}{q}= \frac{1}{r}+\frac{1}{m}$, $a$ is a non-negative number and the dimension $n\ge 1$.
\ehq

\subsection{$L^{q}- L^{q}$ estimates with additional $L^{m}$ regularity for the data} \label{Sec3.3}
From Corollaries \ref{hq4.3.2} and \ref{hq4.3.5} we obtain the following result.
\bmd \label{md4.3.3}
Let $\delta\in \big(\f{\sigma}{2},\sigma \big]$ in (\ref{pt6.3}), $q\in (1,\ity)$ and $m\in [1,q)$. Then the Sobolev solutions to (\ref{pt6.3}) satisfy the following $(L^m \cap L^q)-L^q$ estimates:
\begin{align}
\big\||D|^a u(t,\cdot)\big\|_{L^q}& \lesssim  (1+t)^{(2+[\frac{n}{2}])(1-\frac{\sigma}{2\delta})\frac{1}{r} -\frac{n}{2\delta}(1-\frac{1}{r})-\frac{a}{2\delta}} \|u_0\|_{L^m \cap H^{a,q}} \nonumber\\
& \qquad \quad + (1+t)^{1+(1+[\frac{n}{2}])(1-\frac{\sigma}{2\delta})\frac{1}{r} -\frac{n}{2\delta}(1-\frac{1}{r})-\frac{a}{2\delta}}\|u_1\|_{L^m \cap H^{[a-2\delta]^+,q}}, \label{p4.3.1}\\
\big\||D|^a u_t(t,\cdot)\big\|_{L^q}& \lesssim  (1+t)^{(1+[\frac{n}{2}])(1-\frac{\sigma}{2\delta})\frac{1}{r} -\frac{n}{2\delta}(1-\frac{1}{r})-\frac{a+2(\sigma-\delta)}{2\delta}} \|u_0\|_{L^m \cap H^{a+2(\sigma-\delta),q}} \nonumber\\
& \qquad \quad + (1+t)^{(2+[\frac{n}{2}])(1-\frac{\sigma}{2\delta})\frac{1}{r} -\frac{n}{2\delta}(1-\frac{1}{r})-\frac{a}{2\delta}}\|u_1\|_{L^m \cap H^{a,q}}. \label{p4.3.2}
\end{align}
Moreover, the following $L^q-L^q$ estimates are satisfied:
\begin{align*}
\big\||D|^a u(t,\cdot)\big\|_{L^q}& \lesssim (1+t)^{(2+[\frac{n}{2}])(1-\frac{\sigma}{2\delta})-\frac{a}{2\delta}} \|u_0\|_{H^{a,q}}+ (1+t)^{1+(1+[\frac{n}{2}])(1-\frac{\sigma}{2\delta})-\frac{a}{2\delta}} \|u_1\|_{H^{[a-2\delta]^+,q}}, \\
\big\||D|^a u_t(t,\cdot)\big\|_{L^q}& \lesssim (1+t)^{(1+[\frac{n}{2}])(1-\frac{\sigma}{2\delta})-\frac{a+2(\sigma-\delta)}{2\delta}} \|u_0\|_{H^{a+2(\sigma-\delta),q}}+ (1+t)^{(2+[\frac{n}{2}])(1-\frac{\sigma}{2\delta})-\frac{a}{2\delta}} \|u_1\|_{H^{a,q}}.
\end{align*}
Here $1+ \frac{1}{q}= \frac{1}{r}+\frac{1}{m}$, $a$ is a non-negative number and the dimension $n\ge 1$.
\emd

\begin{nx} \label{nx4.3.3}
\fontshape{n}
\selectfont
The statements in Proposition \ref{md4.3.3} are key tools to prove global (in time) existence results for the semi-linear models (\ref{pt6.1}) and (\ref{pt6.2}). Let us compare the results from Proposition \ref{md4.3.3} to those from Proposition $3.7$ in our previous work \cite{DaoReissig}. First, we can see that there does not appear a singular behavior of the time-dependent coefficients for $t\longrightarrow +0$ in the above estimates. On the one hand, it is important to notice that in (\ref{p4.3.1}) we derived a decay estimate for the fractional derivative of order $a=2\delta$ of the solution with respect to the spatial variables by assuming a suitable higher regularity on $u_0$, that is, $u_0 \in L^m \cap H^{2\delta,q}$, whereas the second data $u_1$ only belongs to $L^m \cap L^q$. This effect does not appear in the case $\delta\in \big(0,\f{\sigma}{2})$. If we assume $u_0 \in L^m \cap H^{a,q}$ for $a>\sigma$, then we choose the second data $u_1$ from the function space  $L^m \cap H^{a-s,q}$. That property brings some benefit in the treatment of the semi-linear models (\ref{pt6.1}) and (\ref{pt6.2}) in the next section. On the other hand, with $a=0$ in (\ref{p4.3.2}) the estimate for the first derivative in time requires less regularity for the data comparing with respect to the estimate for the derivative in space of order $a=\sigma$ in (\ref{p4.3.1}). This property is new in comparison with the case $\delta \in (0,\frac{\sigma}{2})$.
\end{nx}

\section{Treatment of the corresponding semi-linear models} \label{Sec4}

\subsection{Philosophy of our approach} \label{Sec4.1}
In this section, our goal is to apply the estimates from Proposition \ref{md4.3.3} to prove the global (in time) existence of small data Sobolev solutions to the semi-linear models (\ref{pt6.1}) and (\ref{pt6.2}). Some developed tools from Harmonic Analysis come into play such as fractional Gargliardo-Nirenberg inequality from Proposition \ref{fractionalGagliardoNirenberg}, fractional Leibniz rule from Proposition \ref{fractionalLeibniz}, fractional chain rule from Proposition \ref{Propfractionalchainrulegeneral} and fractional Sobolev embedding from Corollary \ref{CorollaryEmbedding}. By recalling the fundamental solutions $K_0$ and $K_1$ defined in Section \ref{Sec3} we write the solutions to (\ref{pt6.3}) in the following form:
$$ u^{ln}(t,x)=K_0(t,x) \ast_{x} u_0(x)+ K_1(t,x) \ast_{x} u_1(x). $$
We apply Duhamel's principle to get the following implicit representation of the solutions to (\ref{pt6.1}) and (\ref{pt6.2}):
$$ u(t,x)= K_0(t,x) \ast_x u_0(x)+ K_1(t,x) \ast_x u_1(x) + \int_0^t K_1(t-\tau,x) \ast_x f(u,u_t) \,d\tau=: u^{ln}(t,x)+ u^{nl}(t,x), $$
where $f(u,u_t)=|u|^p$ or $|u_t|^p$. We choose the data space $\mathcal{A}=\mathcal{A}^{s}_{m,q}$ and introduce the family $\{X(t)\}_{t>0}$ of solution spaces $X(t)$ with the norm
\begin{align*}
\|u\|_{X(t)}:= & \sup_{0\le \tau \le t} \Big( f_{1}(\tau)^{-1}\|u(\tau,\cdot)\|_{L^q}+ f_{\sigma}(\tau)^{-1}\big\||D|^\sigma u(\tau,\cdot)\big\|_{L^q}+ f_{2,s}(\tau)^{-1}\big\||D|^s u(\tau,\cdot)\big\|_{L^q}\\
&\qquad \quad + f_{3}(\tau)^{-1}\|u_t(\tau,\cdot)\|_{L^q}+ f_{4,s}(\tau)^{-1}\big\||D|^{s-2\delta} u_t(\tau,\cdot)\big\|_{L^q} \Big).
\end{align*}
Furthermore, we introduce the family $\{X_0(t)\}_{t>0}$ of space $X_0(t):= C([0,t],H^{s,q})$ with the norm
$$\|w\|_{X_0(t)}:= \sup_{0\le \tau \le t} \Big( f_{1}(\tau)^{-1}\|u(\tau,\cdot)\|_{L^q}+ f_{\sigma}(\tau)^{-1}\big\||D|^\sigma u(\tau,\cdot)\big\|_{L^q}+ f_{2,s}(\tau)^{-1}\big\||D|^s u(\tau,\cdot)\big\|_{L^q} \Big). $$
In both families of spaces we choose the weights
\begin{align*}
&f_{1}(\tau)= (1+\tau)^{1+(1+[\frac{n}{2}])(1-\frac{\sigma}{2\delta})\frac{1}{r}-\frac{n}{2\delta}(1-\frac{1}{r})},\,\, f_{2,s}(\tau)=(1+\tau)^{1+(1+[\frac{n}{2}])(1-\frac{\sigma}{2\delta})\frac{1}{r} -\frac{n}{2\delta}(1-\frac{1}{r})-\frac{s}{2\delta}}, \\
&f_{3}(\tau)=(1+\tau)^{(2+[\frac{n}{2}])(1-\frac{\sigma}{2\delta})\frac{1}{r} -\frac{n}{2\delta}(1-\frac{1}{r})},\,\, f_{4,s}(\tau)=(1+\tau)^{1+(2+[\frac{n}{2}])(1-\frac{\sigma}{2\delta})\frac{1}{r} -\frac{n}{2\delta}(1-\frac{1}{r}) -\frac{s}{2\delta}},
\end{align*}
and
$$ f_{\sigma}(\tau)=(1+\tau)^{1+(1+[\frac{n}{2}])(1-\frac{\sigma}{2\delta})\frac{1}{r} -\frac{n}{2\delta}(1-\frac{1}{r})-\frac{\sigma}{2\delta}}. $$
We define for all $t>0$ the operator $N: \quad u \in X(t) \longrightarrow Nu \in X(t)$ by the formula
$$Nu(t,x)= K_0(t,x) \ast_x u_0(x)+ K_1(t,x) \ast_x u_1(x)+ \int_0^t K_1(t-\tau,x) \ast_x f(u,u_t) \,d\tau. $$
We will prove that the operator $N$ satisfies the following two inequalities:
\begin{align}
\|Nu\|_{X(t)}& \lesssim \|(u_0,u_1)\|_{\mathcal{A}}+ \|u\|^p_{X_0(t)}, \label{pt4.3}\\
\|Nu-Nv\|_{X(t)}& \lesssim \|u-v\|_{X_0(t)} \big(\|u\|^{p-1}_{X_0(t)}+ \|v\|^{p-1}_{X_0(t)}\big). \label{pt4.4}
\end{align}
Then, employing Banach's fixed point theorem leads to local (in time) existence results of large data solutions and global (in time) existence results of small data solutions as well.

\begin{nx}
\fontshape{n}
\selectfont
Replacing $a=s$ and $a=\sigma$ in the statements from Proposition \ref{md4.3.3} and in the definition of the norm of $X(t)$ we conclude
$$ \big\|u^{ln} \big\|_{X(t)} \lesssim \|(u_0,u_1)\|_{\mathcal{A}^{s}_{m,q}} \text{ for }s \ge 0. $$
Hence, in order to prove (\ref{pt4.3}) it is reasonable to indicate the following inequality:
\begin{equation}
\big\|u^{nl} \big\|_{X(t)} \lesssim \|u\|^p_{X_0(t)}. \label{pt4.31}
\end{equation}
\end{nx}

Now we are going to prove our main results from Section \ref{Sec2}.

\subsection{Proof of Theorem \ref{dl2.1}: $s= 2\delta$} \label{Sec4.2}

We introduce the data space $\mathcal{A}:= \mathcal{A}^{2\delta}_{m,q}$ and the solution space
$$X(t):= C([0,t],H^{2\delta,q}) \cap C^1([0,t],L^q), $$
where the weight $f_{4,s}(\tau) \equiv 0$. First, let us prove the inequality (\ref{pt4.31}). In order to control $u^{nl}$, we use the $(L^m \cap L^q)- L^q$ estimates in Proposition \ref{md4.3.3}. Hence, we obtain for $k=0,1$ the following estimates:
$$ \big\||D|^{2\delta k}u^{nl}(t,\cdot)\big\|_{L^q} \lesssim \int_0^{t}(1+t-\tau)^{1+(1+[\frac{n}{2}])(1-\frac{\sigma}{2\delta})\frac{1}{r}-\frac{n}{2\delta}(1-\frac{1}{r})-k}\big\||u(\tau,\cdot)|^p\big\|_{L^m \cap L^q} \,d\tau. $$
Therefore, it is necessary to estimate $|u(\tau,x)|^p$ in $L^m \cap L^q$. We proceed as follows:
$$\big\||u(\tau,\cdot)|^p\big\|_{L^m \cap L^q} \lesssim \|u(\tau,\cdot)\|^p_{L^{mp}}+ \|u(\tau,\cdot)\|^p_{L^{qp}}. $$
Employing the fractional Gagliardo-Nirenberg inequality from Proposition \ref{fractionalGagliardoNirenberg} leads to
$$\big\||u(\tau,\cdot)|^p\big\|_{L^m \cap L^q} \lesssim (1+\tau)^{p\big(1+(1+[\frac{n}{2}])(1-\frac{\sigma}{2\delta})\frac{1}{r}-\frac{n}{2\delta}(\frac{1}{m}-\frac{1}{mp})\big)}\|u\|^p_{X_0(\tau)} $$
provided that the condition (\ref{GN2A}) is fulfilled. From the above estimate we get
$$\big\||D|^{2\delta k}u^{nl}(t,\cdot)\big\|_{L^q} \lesssim \|u\|^p_{X_0(t)} \int_0^{t}(1+t-\tau)^{1+(1+[\frac{n}{2}])(1-\frac{\sigma}{2\delta})\frac{1}{r}-\frac{n}{2\delta}(1-\frac{1}{r})-k} (1+\tau)^{p\big(1+(1+[\frac{n}{2}])(1-\frac{\sigma}{2\delta})\frac{1}{r}-\frac{n}{2\delta}(\frac{1}{m}-\frac{1}{mp})\big)} d\tau. $$
The key tool relies now in Lemma \ref{LemmaIntegral}. Because of the condition (\ref{exponent2A}), applying Lemma \ref{LemmaIntegral} by choosing
$\alpha= -1-(1+[\frac{n}{2}])(1-\frac{\sigma}{2\delta})\frac{1}{r}+\frac{n}{2\delta}(1-\frac{1}{r})+k$ and
$\beta= p\big(-1-(1+[\frac{n}{2}])(1-\frac{\sigma}{2\delta})\frac{1}{r}+\frac{n}{2\delta}(\frac{1}{m}-\frac{1}{mp})\big)$ we derive
\begin{align*}
&\int_0^{t}(1+t-\tau)^{1+(1+[\frac{n}{2}])(1-\frac{\sigma}{2\delta})\frac{1}{r}-\frac{n}{2\delta}(1-\frac{1}{r})-k} (1+\tau)^{p\big(1+(1+[\frac{n}{2}])(1-\frac{\sigma}{2\delta})\frac{1}{r}-\frac{n}{2\delta}(\frac{1}{m}-\frac{1}{mp})\big)} d\tau\\
&\qquad \quad \lesssim (1+t)^{1+(1+[\frac{n}{2}])(1-\frac{\sigma}{2\delta})\frac{1}{r}-\frac{n}{2\delta}(1-\frac{1}{r})-k}.
\end{align*}
As a result, we may conclude for $k=0,1$ the estimates
\begin{equation}
\big\||D|^{2\delta k}u^{nl}(t,\cdot)\big\|_{L^q} \lesssim (1+t)^{1+(1+[\frac{n}{2}])(1-\frac{\sigma}{2\delta})\frac{1}{r}-\frac{n}{2\delta}(1-\frac{1}{r})-k} \|u\|^p_{X_0(t)}. \label{t2.1.1}
\end{equation}
Analogously, we also obtain
\begin{align}
\big\||D|^{\sigma}u^{nl}(t,\cdot)\big\|_{L^q} &\lesssim (1+t)^{1+(1+[\frac{n}{2}])(1-\frac{\sigma}{2\delta})\frac{1}{r}-\frac{n}{2\delta}(1-\frac{1}{r})-\frac{\sigma}{2\delta}} \|u\|^p_{X_0(t)}, \label{t2.1.2}\\
\big\|\partial_t u^{nl}(t,\cdot)\big\|_{L^q} &\lesssim (1+t)^{(2+[\frac{n}{2}])(1-\frac{\sigma}{2\delta})\frac{1}{r} -\frac{n}{2\delta}(1-\frac{1}{r})}\|u\|^p_{X_0(t)}. \label{t2.1.3}
\end{align}
From (\ref{t2.1.1}) to (\ref{t2.1.3}) and the definition of the norm in $X(t)$, we arrive immediately at the inequality (\ref{pt4.31}). \medskip

\noindent Next, let us prove the estimate (\ref{pt4.4}). Using again the $(L^m \cap L^q)- L^q$ estimates in Proposition \ref{md4.3.3}, we have for two functions $u$ and $v$ from $X(t)$ the following estimate for $k=0,1$:
$$ \big\||D|^{2\delta k}\big(Nu(t,\cdot)-Nv(t,\cdot)\big)\big\|_{L^q} \lesssim \int_0^{t}(1+t-\tau)^{1+(1+[\frac{n}{2}])(1-\frac{\sigma}{2\delta})\frac{1}{r}-\frac{n}{2\delta}(1-\frac{1}{r})-k}\big\||u(\tau,\cdot)|^p-|v(\tau,\cdot)|^p\big\|_{L^m \cap L^q}\,d\tau. $$
Applying H\"{o}lder's inequality leads to
\begin{align*}
\big\||u(\tau,\cdot)|^p-|v(\tau,\cdot)|^p\big\|_{L^q}& \lesssim \|u(\tau,\cdot)-v(\tau,\cdot)\|_{L^{qp}} \big(\|u(\tau,\cdot)\|^{p-1}_{L^{qp}}+\|v(\tau,\cdot)\|^{p-1}_{L^{qp}}\big),\\
\big\||u(\tau,\cdot)|^p-|v(\tau,\cdot)|^p\big\|_{L^m}& \lesssim \|u(\tau,\cdot)-v(\tau,\cdot)\|_{L^{mp}} \big(\|u(\tau,\cdot)\|^{p-1}_{L^{mp}}+\|v(\tau,\cdot)\|^{p-1}_{L^{mp}}\big).
\end{align*}
In the same way as the proof of (\ref{pt4.3}), employing the fractional Gagliardo-Nirenberg inequality from Proposition \ref{fractionalGagliardoNirenberg} to the terms
$$ \|u(\tau,\cdot)-v(\tau,\cdot)\|_{L^\eta }, \text{ }\|u(\tau,\cdot)\|_{L^\eta}, \text{ }\|v(\tau,\cdot)\|_{L^\eta} $$
with $\eta=qp$ and $\eta=mp$ we have for $k=0,1$ the estimates
$$\big\||D|^{2\delta k}\big(Nu(t,\cdot)-Nv(t,\cdot)\big)\big\|_{L^q} \lesssim (1+t)^{1+(1+[\frac{n}{2}])(1-\frac{\sigma}{2\delta})\frac{1}{r}-\frac{n}{2\delta}(1-\frac{1}{r})-k}\|u-v\|_{X_0(t)}
\big(\|u\|^{p-1}_{X_0(t)}+\|v\|^{p-1}_{X_0(t)}\big). $$
Analogously, we also derive
\begin{align*}
\big\||D|^{\sigma}\big(Nu(t,\cdot)-Nv(t,\cdot)\big)\big\|_{L^q} &\lesssim (1+t)^{1+(1+[\frac{n}{2}])(1-\frac{\sigma}{2\delta})\frac{1}{r}-\frac{n}{2\delta}(1-\frac{1}{r})-\frac{\sigma}{2\delta}}\|u-v\|_{X_0(t)}\big(\|u\|^{p-1}_{X_0(t)}+\|v\|^{p-1}_{X_0(t)}\big), \\
\big\|\partial_t\big(Nu(t,\cdot)-Nv(t,\cdot)\big)\big\|_{L^q} &\lesssim (1+t)^{(2+[\frac{n}{2}])(1-\frac{\sigma}{2\delta})\frac{1}{r} -\frac{n}{2\delta}(1-\frac{1}{r})}\|u-v\|_{X_0(t)}\big(\|u\|^{p-1}_{X_0(t)}+\|v\|^{p-1}_{X_0(t)}\big).
\end{align*}
From the definition of the norm in $X(t)$, we may conclude the inequality (\ref{pt4.4}). Summarizing, the proof to Theorem \ref{dl2.1} is complete.

\subsection{Proof of Theorem \ref{dl2.2}: $0 < s <2\delta$} \label{Sec4.3}

We introduce the data space $\mathcal{A}:= \mathcal{A}^{s}_{m,q}$ and the solution space
$$X(t):= C([0,t],H^{s,q}), $$
where the weights $f_{\sigma}(\tau)=f_{3}(\tau)=f_{4,s}(\tau) \equiv 0$. We can notice that $X_0(t)$ and $X(t)$ coincide in (\ref{pt4.4}) and (\ref{pt4.31}). In order to prove these two inequalities, we use the $(L^m \cap L^q)- L^q$ estimates from Proposition \ref{md4.3.3}. Hence, we derive for $k=0,1$ the following estimates:
\begin{align*}
 \big\||D|^{ks}u^{nl}(t,\cdot)\big\|_{L^q} &\lesssim \int_0^{t}(1+t-\tau)^{1+(1+[\frac{n}{2}])(1-\frac{\sigma}{2\delta})\frac{1}{r}-\frac{n}{2\delta}(1-\frac{1}{r})-\frac{ks}{2\delta}}
 \big\||u(\tau,\cdot)|^p\big\|_{L^m \cap L^q} \,d\tau, \\
\big\||D|^{ks}\big(Nu(t,\cdot)-Nv(t,\cdot)\big)\big\|_{L^q} &\lesssim \int_0^{t}(1+t-\tau)^{1+(1+[\frac{n}{2}])(1-\frac{\sigma}{2\delta})\frac{1}{r}-\frac{n}{2\delta}(1-\frac{1}{r})-\frac{ks}{2\delta}}
\big\||u(\tau,\cdot)|^p-|v(\tau,\cdot)|^p\big\|_{L^m \cap L^q}\,d\tau.
\end{align*}
In the same way as we did in the proof of Theorem \ref{dl2.1} we obtain for $k=0,1$ the following estimates:
\begin{align*}
\big\||D|^{ks}u^{nl}(t,\cdot)\big\|_{L^q} &\lesssim (1+t)^{1+(1+[\frac{n}{2}])(1-\frac{\sigma}{2\delta})\frac{1}{r}-\frac{n}{2\delta}(1-\frac{1}{r})-\frac{ks}{2\delta}}\|u\|^p_{X(t)},\\
\big\||D|^{ks}\big(Nu(t,\cdot)-Nv(t,\cdot)\big)\big\|_{L^q} &\lesssim (1+t)^{1+(1+[\frac{n}{2}])(1-\frac{\sigma}{2\delta})\frac{1}{r}-\frac{n}{2\delta}(1-\frac{1}{r})-\frac{ks}{2\delta}}\|u-v\|_{X(t)}\big(\|u\|^{p-1}_{X(t)}+\|v\|^{p-1}_{X(t)}\big),
\end{align*}
provided that the conditions (\ref{exponent3A}) and (\ref{GN3A}) are fulfilled. From the definition of the norm in $X(t)$ we can conclude immediately the inequalities (\ref{pt4.31}) and (\ref{pt4.4}). Summarizing, Theorem \ref{dl2.2} is proved.

\subsection{Proof of Theorem \ref{dl2.3}: $2\delta< s \le 2\delta + \frac{n}{q}$} \label{Sec4.4}

We introduce the data space $\mathcal{A}:= \mathcal{A}^{s}_{m,q}$ and the solution space
$$X(t):= C([0,t],H^{s,q}) \cap C^1([0,t],H^{s-2\delta,q}), $$
where the weight $f_{\sigma}(\tau) \equiv 0$. First, let us prove the inequality (\ref{pt4.31}). We need to control all norms
$$\|u^{nl}(t,\cdot)\|_{L^q},\,\, \|u_t^{nl}(t,\cdot)\|_{L^q},\,\, \big\||D|^{s}u^{nl}(t,\cdot)\big\|_{L^q},\,\, \big\||D|^{s-2\delta}u_t^{nl}(t,\cdot)\big\|_{L^q}. $$
In the analogous way as we did in the proof of Theorem \ref{dl2.1}, we derive the following estimates:
\begin{align}
\|u^{nl}(t,\cdot)\|_{L^q} &\lesssim (1+t)^{1+(1+[\frac{n}{2}])(1-\frac{\sigma}{2\delta})\frac{1}{r}-\frac{n}{2\delta}(1-\frac{1}{r})} \|u\|^p_{X_0(t)}, \label{t4A1}\\
\|\partial_t u^{nl}(t,\cdot)\|_{L^q} &\lesssim  (1+\tau)^{(2+[\frac{n}{2}])(1-\frac{\sigma}{2\delta})\frac{1}{r}-\frac{n}{2\delta}(1-\frac{1}{r})} \|u\|^p_{X_0(t)}, \label{t4A2}
\end{align}
provided that the condition ({\ref{exponent4A}}) is satisfied and
\begin{equation}
p\in \Big[\frac{q}{m},\ity \Big)  \text{ if } n\le qs, \text{ or }p \in \Big[\f{q}{m}, \f{n}{n-qs}\Big] \text{ if } n \in \Big(qs, \f{q^2 s}{q-m}\Big]. \label{t4A3}
\end{equation}

\noindent Now, let us turn to estimate the norm $\big\||D|^{s-2\delta} u_t^{nl}(t,\cdot)\big\|_{L^q}$. We use the $(L^m \cap L^q)- L^q$ estimates from Proposition \ref{md4.3.3} to get
$$ \big\||D|^{s-2\delta} u_t^{nl}(t,\cdot)\big\|_{L^q} \lesssim \int_0^{t} (1+t-\tau)^{1+(2+[\frac{n}{2}])(1-\frac{\sigma}{2\delta})\frac{1}{r} -\frac{n}{2\delta}(1-\frac{1}{r}) -\frac{s}{2\delta}}\big\||u(\tau,\cdot)|^p\big\|_{L^m \cap L^q \cap \dot{H}^{s-2\delta, q}} \,d\tau. $$
The integrals with $\big\||u(\tau,\cdot)|^p\big\|_{L^m \cap L^q}$ and $\big\||u(\tau,\cdot)|^p\big\|_{L^ q}$ will be handled as we did to obtain (\ref{t4A1}). In order to estimate $\big\||u(\tau,\cdot)|^p\big\|_{\dot{H}^{s-2\delta, q}}$, we shall apply the fractional chain rule with $p> \lceil s-2\delta \rceil$ from Proposition \ref{Propfractionalchainrulegeneral} and the fractional Gagliardo-Nirenberg inequality from Proposition \ref{fractionalGagliardoNirenberg}. Hence, we get
\begin{align*}
&\big\||u(\tau,\cdot)|^p\big\|_{\dot{H}^{s-2\delta, q}} \lesssim \|u(\tau,\cdot)\|^{p-1}_{L^{q_1}}\,\,\big\||D|^{s-2\delta}u(\tau,\cdot)\big\|_{L^{q_2}}\\
&\qquad \lesssim \|u(\tau,\cdot)\|^{(p-1)(1-\theta_{q_1})}_{L^q}\,\,\big\||D|^s u(\tau,\cdot)\big\|^{(p-1)\theta_{q_1}}_{L^q}\,\,\|u(\tau,\cdot)\|^{1-\theta_{q_2}}_{L^q}\,\,\big\||D|^s u(\tau,\cdot)\big\|^{\theta_{q_2}}_{L^q}\\
&\qquad \lesssim (1+\tau)^{p(1+(1+[\frac{n}{2}])(1-\frac{\sigma}{2\delta})\frac{1}{r}-\frac{n}{2\delta}(\frac{1}{m}-\frac{1}{qp}))-\frac{s-2\delta}{2\delta}}\|u\|^p_{X_0(\tau)},
\end{align*}
where
$$\frac{p-1}{q_1}+\frac{1}{q_2}= \frac{1}{q},\,\, \theta_{q_1}= \frac{n}{s}\Big(\frac{1}{q}-\frac{1}{q_1}\Big) \in [0,1],\,\, \theta_{q_2}= \frac{n}{s}\Big(\frac{1}{q}-\frac{1}{q_2}+\frac{s-2\delta}{n}\Big) \in \Big[\frac{s-2\delta}{s},1\Big]. $$
These conditions imply the restriction
\begin{equation}
1<p\le 1+\frac{q2\delta}{n-qs} \text{ if } n>qs, \text{ or } p>1 \text{ if } n \le qs. \label{t4A4}
\end{equation}
Therefore, we obtain
\begin{equation}
\big\||D|^{s-2\delta} u_t^{nl}(t,\cdot)\big\|_{L^q} \lesssim (1+t)^{1+(2+[\frac{n}{2}])(1-\frac{\sigma}{2\delta})\frac{1}{r} -\frac{n}{2\delta}(1-\frac{1}{r}) -\frac{s}{2\delta}} \|u\|^p_{X_0(t)}. \label{t4A5}
\end{equation}
In the analogous way we also derive
\begin{equation}
\big\||D|^s u^{nl}(t,\cdot)\big\|_{L^q} \lesssim (1+t)^{1+(1+[\frac{n}{2}])(1-\frac{\sigma}{2\delta})\frac{1}{r} -\frac{n}{2\delta}(1-\frac{1}{r})-\frac{s}{2\delta}} \|u\|^p_{X_0(t)}. \label{t4A6}
\end{equation}
Summarizing, from (\ref{t4A1}) to (\ref{t4A2}), (\ref{t4A5}) to (\ref{t4A6}) and the definition of the norm in $X(t)$ the inequality (\ref{pt4.31}) follows immediately.\medskip

\noindent Next, let us prove the inequality (\ref{pt4.4}). Following the proof of Theorem \ref{dl2.1}, the new difficulty is to control the term $\big\||u(\tau,\cdot)|^p-|v(\tau,\cdot)|^p\big\|_{\dot{H}^{s-2\delta,q}}$. The integral representation
$$ |u(\tau,x)|^p-|v(\tau,x)|^p=p\int_0^1 \big(u(\tau,x)-v(\tau,x)\big)G\big(\omega u(\tau,x)+(1-\omega)v(\tau,x)\big)\,d\omega, $$
where $G(u)=u|u|^{p-2}$, leads to
$$\big\||u(\tau,\cdot)|^p-|v(\tau,\cdot)|^p\big\|_{\dot{H}^{s-2\delta,q}} \lesssim \int_0^1 \Big\||D|^{s-2\delta}\Big(\big(u(\tau,\cdot)-v(\tau,\cdot)\big)G\big(\omega u(\tau,\cdot)+(1-\omega)v(\tau,\cdot)\big)\Big)\Big\|_{L^q}\,d\omega. $$
Applying the fractional Leibniz formula from Proposition \ref{fractionalLeibniz} we derive the following estimate:
\begin{align*}
& \big\||u(\tau,\cdot)|^p-|v(\tau,\cdot)|^p\big\|_{\dot{H}^{s-2\delta,q}} \\ & \qquad \lesssim \big\||D|^{s-2\delta}\big(u(\tau,\cdot)-v(\tau,\cdot)\big)\big\|_{L^{r_1}} \int_0^1 \big\|G\big(\omega u(\tau,\cdot)+(1-\omega)v(\tau,\cdot)\big)\big\|_{L^{r_2}}\,d\omega\\
& \quad \qquad + \|u(\tau,\cdot)-v(\tau,\cdot)\|_{L^{r_3}}  \int_0^1 \big\||D|^{s-2\delta}G\big(\omega u(\tau,\cdot)+(1-\omega)v(\tau,\cdot)\big)\big\|_{L^{r_4}}\,d\omega\\
& \qquad \lesssim \big\||D|^{s-2\delta}\big(u(\tau,\cdot)-v(\tau,\cdot)\big)\big\|_{L^{r_1}} \Big(\|u(\tau,\cdot)\|^{p-1}_{L^{r_2 (p-1)}}+ \|v(\tau,\cdot)\|^{p-1}_{L^{r_2 (p-1)}}\Big)\\
&\qquad \quad + \|u(\tau,\cdot)-v(\tau,\cdot)\|_{L^{r_3}}  \int_0^1 \big\||D|^{s-2\delta}G\big(\omega u(\tau,\cdot)+(1-\omega)v(\tau,\cdot)\big)\big\|_{L^{r_4}}\,d\omega,
\end{align*}
where
$$\frac{1}{r_1}+\frac{1}{r_2}= \frac{1}{r_3}+\frac{1}{r_4}= \frac{1}{q}.$$
Taking into consideration the fractional Gargliardo- Nirenberg inequality from Proposition \ref{fractionalGagliardoNirenberg} we obtain
\begin{align*}
\big\||D|^{s-2\delta}\big(u(\tau,\cdot)-v(\tau,\cdot)\big)\big\|_{L^{r_1}}&\lesssim \|u(\tau,\cdot)-v(\tau,\cdot)\|^{\theta_1}_{\dot{H}^{s,q}}\,\,\|u(\tau,\cdot)-v(\tau,\cdot)\|^{1-\theta_1}_{L^q},\\
\|u(\tau,\cdot)\|_{L^{r_2 (p-1)}}&\lesssim \|u(\tau,\cdot)\|^{\theta_2}_{\dot{H}^{s,q}}\,\,\|u(\tau,\cdot)\|^{1-\theta_2}_{L^q},\\
\|u(\tau,\cdot)-v(\tau,\cdot)\|_{L^{r_3}}&\lesssim \|u(\tau,\cdot)-v(\tau,\cdot)\|^{\theta_3}_{\dot{H}^{s,q}}\,\,\|u(\tau,\cdot)-v(\tau,\cdot)\|^{1-\theta_3}_{L^q},
\end{align*}
where
$$\theta_1= \frac{n}{s}\Big(\frac{1}{q}-\frac{1}{r_1}+\frac{s-2\delta}{n}\Big) \in \Big[\frac{s-2\delta}{s},1\Big],\,\, \theta_2= \frac{n}{s}\Big(\frac{1}{q}-\frac{1}{r_2(p-1)}\Big) \in [0,1],\,\, \theta_3= \frac{n}{s}\Big(\frac{1}{q}-\frac{1}{r_3}\Big) \in [0,1]. $$
Because $\omega \in [0,1]$ is a parameter, employing again the fractional chain rule with $p >1+ \lceil s-2\delta \rceil$ from Proposition \ref{Propfractionalchainrulegeneral} and the fractional Gagliardo- Nirenberg inequality we get
\begin{align*}
&\big\||D|^{s-2\delta}G\big(\omega u(\tau,\cdot)+(1-\omega)v(\tau,\cdot)\big)\big\|_{L^{r_4}}\\
&\qquad \lesssim \|\omega u(\tau,\cdot)+(1-\omega)v(\tau,\cdot)\|^{p-2}_{L^{r_5}}\,\, \big\||D|^{s-2\delta}\big(\omega u(\tau,\cdot)+(1-\omega)v(\tau,\cdot)\big)\big\|_{L^{r_6}}\\
&\qquad \lesssim \|\omega u(\tau,\cdot)+(1-\omega)v(\tau,\cdot)\|^{(p-2)\theta_5+\theta_6}_{\dot{H}^{s,q}}\,\, \|\omega u(\tau,\cdot)+(1-\omega)v(\tau,\cdot)\|^{(p-2)(1-\theta_5)+1-\theta_6}_{L^q},
\end{align*}
where
$$\frac{p-2}{r_5}+\frac{1}{r_6}= \frac{1}{r_4},\,\, \theta_5= \frac{n}{s}\Big(\frac{1}{q}-\frac{1}{r_5}\Big) \in [0,1],\,\, \theta_6= \frac{n}{s}\Big(\frac{1}{q}-\frac{1}{r_6}+\frac{s-2\delta}{n}\Big) \in \Big[\frac{s-2\delta}{s},1\Big]. $$
All together it follows
\begin{align*}
&\int_0^1 \big\||D|^{s-2\delta}G\big(\omega u(\tau,\cdot)+(1-\omega)v(\tau,\cdot)\big)\big\|_{L^{r_4}}d\omega\\
&\quad \lesssim \big(\|u(\tau,\cdot)\|_{\dot{H}^{s,q}}+\|v(\tau,\cdot)\|_{\dot{H}^{s,q}}\big)^{(p-2)\theta_5+\theta_6} \big(\|u(\tau,\cdot)\|_{L^q}+\|v(\tau,\cdot)\|_{L^q} \big)^{(p-2)(1-\theta_5)+1-\theta_6}.
\end{align*}
Hence, we derived the following estimate:
\begin{eqnarray*} && \big\||u(\tau,\cdot)|^p-|v(\tau,\cdot)|^p\big\|_{\dot{H}^{s-2\delta,q}} \\ && \qquad \lesssim (1+\tau)^{p\big(1+(1+[\frac{n}{2}])(1-\frac{\sigma}{2\delta})\frac{1}{r}-\frac{n}{2\delta}(\frac{1}{m}-\frac{1}{qp})\big)-\frac{s-2\delta}{2\delta}}\|u-v\|_{X_0(\tau)}\big( \|u\|^{p-1}_{X_0(\tau)}+ \|v\|^{p-1}_{X_0(\tau)} \big),\end{eqnarray*}
where we note that
$$\theta_1+ (p-1)\theta_2= \theta_3+ (p-2)\theta_5+ \theta_6= \frac{n}{s}\Big(\frac{p-1}{q}+\frac{s-2\delta}{n}\Big). $$
Therefore, we have proved that
\begin{align*}
\big\||D|^{s-2\delta} \partial_t \big(Nu(t,\cdot)- Nv(t,\cdot)\big)\big\|_{L^q} &\lesssim  (1+t)^{1+(2+[\frac{n}{2}])(1-\frac{\sigma}{2\delta})\frac{1}{r} -\frac{n}{2\delta}(1-\frac{1}{r}) -\frac{s}{2\delta}} \|u-v\|_{X_0(t)}\big( \|u\|^{p-1}_{X_0(t)}+ \|v\|^{p-1}_{X_0(t)} \big), \\
\big\||D|^s \big(Nu(t,\cdot)- Nv(t,\cdot)\big)\big\|_{L^q} &\lesssim (1+t)^{1+(1+[\frac{n}{2}])(1-\frac{\sigma}{2\delta})\frac{1}{r} -\frac{n}{2\delta}(1-\frac{1}{r})-\frac{s}{2\delta}} \|u-v\|_{X_0(t)}\big( \|u\|^{p-1}_{X_0(t)}+ \|v\|^{p-1}_{X_0(t)} \big).
\end{align*}
From the definition of the norm in $X(t)$ the inequality (\ref{pt4.4}) follows. Summarizing, the proof of Theorem \ref{dl2.3} is complete.

\begin{nx}
\fontshape{n}
\selectfont
It is clear that one should explain if one can really choose parameters $q_1$, $q_2$, $r_1,\cdots, r_6$ and $\theta_1,\cdots,\theta_6$ as required in the proof to Theorem \ref{dl2.3}. Following the explanations as we did in Remark $4.2$ in \cite{DaoReissig} we may conclude the conditions
\begin{equation}
2\le p \le 1+\frac{q2\delta}{n-qs} \text{ if } n>qs, \text{ or }p \ge 2 \text{ if }n \le qs. \label{Conditionequation}
\end{equation}
These conditions are sufficient to guarantee the existence of all these parameters satisfying the required conditions.
\end{nx}

\subsection{Proof of Theorem \ref{dl2.4}: $s>2\delta+\frac{n}{q}$} \label{Sec4.5}

We introduce both spaces of the data and the solutions as in Theorem \ref{dl2.3}, where the weight $f_{\sigma}(\tau) \equiv 0$. On the one hand, we can repeat exactly how to estimate the terms $|u(\tau,\cdot)|^p$ and $|u(\tau,\cdot)|^p-|v(\tau,\cdot)|^p$ in $L^m$ and $L^q$ as we did in the proof to Theorem \ref{dl2.3}. On the other hand, let us control the above terms in $\dot{H}^{s-2\delta, q}$ by applying the fractional powers rule and the fractional Sobolev embedding.\medskip

First, let us begin to estimate $\big\||u(\tau,\cdot)|^p\big\|_{\dot{H}^{s-2\delta, q}}$. We apply Corollary \ref{Corfractionalhomogeneous} for fractional powers with $s-2\delta \in \big(\frac{n}{q},p\big)$ and Corollary \ref{CorollaryEmbedding} with a suitable $s^* <\frac{n}{q}$ to derive
$$ \big\||u(\tau,\cdot)|^p\big\|_{\dot{H}^{s-2\delta, q}}\lesssim \|u(\tau,\cdot)\|_{\dot{H}^{s-2\delta, q}}\|u(\tau,\cdot)\|^{p-1}_{L^\ity} \lesssim \|u(\tau,\cdot)\|_{\dot{H}^{s-2\delta, q}}\big(\|u(\tau,\cdot)\|_{\dot{H}^{s^*, q}}+ \|u(\tau,\cdot)\|_{\dot{H}^{s-2\delta, q}}\big)^{p-1}. $$
Using again the fractional Gagliardo-Nirenberg inequality from Proposition \ref{fractionalGagliardoNirenberg} gives
\begin{align*}
\|u(\tau,\cdot)\|_{\dot{H}^{s-2\delta, q}} &\lesssim \|u(\tau,\cdot)\|^{1-\theta_1}_{L^q}\,\, \big\||D|^s u(\tau,\cdot)\big\|^{\theta_1}_{L^q} \lesssim (1+\tau)^{1+(1+[\frac{n}{2}])(1-\frac{\sigma}{2\delta})\frac{1}{r}-\frac{n}{2\delta}(1-\frac{1}{r})-\frac{s-2\delta}{2\delta}}\|u\|_{X_0(\tau)}, \\
\|u(\tau,\cdot)\|_{\dot{H}^{s^*, q}} &\lesssim \|u(\tau,\cdot)\|^{1-\theta_2}_{L^q}\,\, \big\||D|^s u(\tau,\cdot)\big\|^{\theta_2}_{L^q} \lesssim (1+\tau)^{1+(1+[\frac{n}{2}])(1-\frac{\sigma}{2\delta})\frac{1}{r}-\frac{n}{2\delta}(1-\frac{1}{r})-\frac{s^*}{2\delta}}\|u\|_{X_0(\tau)},
\end{align*}
where $\theta_1= 1- \frac{2\delta}{s}$ and $\theta_2= \frac{s^*}{s}$. Therefore, we obtain
\begin{align*}
&\big\||u(\tau,\cdot)|^p\big\|_{\dot{H}^{s-2\delta, q}}\lesssim (1+\tau)^{p(1+(1+[\frac{n}{2}])(1-\frac{\sigma}{2\delta})\frac{1}{r} -\frac{n}{2\delta}(1-\frac{1}{r}))-\frac{s-2\delta}{2\delta}- (p-1)\frac{s^*}{2\delta}} \|u\|_{X_0(\tau)} \\
&\qquad \lesssim (1+\tau)^{p(1+(1+[\frac{n}{2}])(1-\frac{\sigma}{2\delta})\frac{1}{r}-\frac{n}{2\delta}(\frac{1}{m}-\frac{1}{mp}))}\|u\|^p_{X_0(\tau)}, \end{align*}
if we choose  $s^*= \frac{n}{q}-\e$ with a sufficiently small positive $\e$.\\
Next, let us turn to estimate the term $\big\||u(\tau,\cdot)|^p-|v(\tau,\cdot)|^p\big\|_{\dot{H}^{s-2\delta,q}}$. Then, repeating the corresponding steps of the proof of Theorem \ref{dl2.3} and using analogous arguments as in the first step we conclude
$$\big\||u(\tau,\cdot)|^p-|v(\tau,\cdot)|^p\big\|_{\dot{H}^{s-2\delta,q}} \lesssim (1+\tau)^{p(1+(1+[\frac{n}{2}])(1-\frac{\sigma}{2\delta})\frac{1}{r}-\frac{n}{2\delta}(\frac{1}{m}-\frac{1}{mp}))} \|u-v\|_{X_0(t)}\big( \|u\|^{p-1}_{X_0(t)}+ \|v\|^{p-1}_{X_0(t)} \big), $$
provided that the conditions $p>2$ and $p> 1+s-2\delta$ are satisfied. Summarizing, Theorem \ref{dl2.4} is proved.

\subsection{Proof of Theorem \ref{dl2.5}: $s>2\delta+\frac{n}{q}$} \label{Sec4.6}

We introduce both spaces for the data and the solutions as in Theorem \ref{dl2.3}, where the weight $f_{\sigma}(\tau) \equiv 0$. But now the space $X_0(t)$ is replaced by the space $X(t)$ in both inequalities (\ref{pt4.4}) and (\ref{pt4.31}). First, let us prove the inequality (\ref{pt4.31}). In order to estimate $u^{nl}$, we apply the $L^m \cap L^q- L^q$ estimates from Proposition \ref{md4.3.3}. Therefore, we derive
$$ \big\|u^{nl}(t,\cdot)\big\|_{L^q} \lesssim \int_0^{t}(1+t-\tau)^{1+(1+[\frac{n}{2}])(1-\frac{\sigma}{2\delta})\frac{1}{r}-\frac{n}{2\delta}(1-\frac{1}{r})}\big\||u_t(\tau,\cdot)|^p\big\|_{L^m \cap L^q}\, d\tau. $$
Moreover, we get
$$ \big\||u_t(\tau,\cdot)|^p\big\|_{L^m \cap L^q} \lesssim \|u_t(\tau,\cdot)\|^p_{L^{mp}}+ \|u_t(\tau,\cdot)\|^p_{L^{qp}}. $$
After applying the fractional Gagliardo-Nirenberg inequality from Proposition \ref{fractionalGagliardoNirenberg}
we arrive at
$$ \big\||u_t(\tau,\cdot)|^p\big\|_{L^m \cap L^q} \lesssim (1+\tau)^{p((2+[\frac{n}{2}])(1-\frac{\sigma}{2\delta})\frac{1}{r}-\frac{n}{2\delta}(\frac{1}{m}-\frac{1}{mp}))}\|u\|^p_{X(\tau)}, $$
provided that the condition $p \in \big[\frac{q}{m},\ity \big)$ is fulfilled due to $s > 2\delta+\frac{n}{q}$. Hence, we obtain
$$ \|u^{nl}(t,\cdot)\|_{L^q} \lesssim \|u\|^p_{X(t)} \int_0^{t}(1+t-\tau)^{1+(1+[\frac{n}{2}])(1-\frac{\sigma}{2\delta})\frac{1}{r}-\frac{n}{2\delta}(1-\frac{1}{r})} (1+\tau)^{p((2+[\frac{n}{2}])(1-\frac{\sigma}{2\delta})\frac{1}{r}-\frac{n}{2\delta}(\frac{1}{m}-\frac{1}{mp}))}\, d\tau. $$
The key tool relies now in using Lemma \ref{LemmaIntegral}. Because of condition (\ref{exponent6A}), applying Lemma \ref{LemmaIntegral} after choosing
\[ \alpha= -1-\Big(1+\Big[\frac{n}{2}\Big]\Big)\Big(1-\frac{\sigma}{2\delta}\Big)\frac{1}{r}+\frac{n}{2\delta}\Big(1-\frac{1}{r}\Big)\,\,\,\mbox{ and}\,\,\,
\beta= p\Big(-\Big(2+\Big[\frac{n}{2}\Big]\Big)\Big(1-\frac{\sigma}{2\delta}\Big)\frac{1}{r}+\frac{n}{2\delta}\Big(\frac{1}{m}-\frac{1}{mp}\Big)\Big),\] we have
\begin{eqnarray*} && \int_0^{t}(1+t-\tau)^{1+(1+[\frac{n}{2}])(1-\frac{\sigma}{2\delta})\frac{1}{r}-\frac{n}{2\delta}(1-\frac{1}{r})} (1+\tau)^{p((2+[\frac{n}{2}])(1-\frac{\sigma}{2\delta})\frac{1}{r}-\frac{n}{2\delta}(\frac{1}{m}-\frac{1}{mp}))} \, d\tau \\ && \qquad \lesssim (1+t)^{1+(1+[\frac{n}{2}])(1-\frac{\sigma}{2\delta})\frac{1}{r}-\frac{n}{2\delta}(1-\frac{1}{r})}.\end{eqnarray*}
As a result, we arrive at the following estimate:
\begin{equation}
\|u^{nl}(t,\cdot)\|_{L^q} \lesssim (1+t)^{1+(1+[\frac{n}{2}])(1-\frac{\sigma}{2\delta})\frac{1}{r}-\frac{n}{2\delta}(1-\frac{1}{r})} \|u\|^p_{X(t)}. \label{t6A1}
\end{equation}
Analogously, we also derive
\begin{equation}
\|\partial_t u^{nl}(t,\cdot)\|_{L^q}\lesssim (1+\tau)^{(2+[\frac{n}{2}])(1-\frac{\sigma}{2\delta})\frac{1}{r}-\frac{n}{2\delta}(1-\frac{1}{r})} \|u\|^p_{X(t)}. \label{t6A2}
\end{equation}
Now, let us control the norm $\big\||D|^{s-2\delta} u_t^{nl}(t,\cdot)\big\|_{L^q}$. We get
$$ \big\||D|^{s-2\delta} u_t^{nl}(t,\cdot)\big\|_{L^q} \lesssim \int_0^{t} (1+t-\tau)^{1+(2+[\frac{n}{2}])(1-\frac{\sigma}{2\delta})\frac{1}{r} -\frac{n}{2\delta}(1-\frac{1}{r}) -\frac{s}{2\delta}}\big\||u_t(\tau,\cdot)|^p\big\|_{L^m \cap L^q \cap \dot{H}^{s-2\delta, q}} d\tau. $$
The integrals with estimates for $\big\||u_t(\tau,\cdot)|^p\big\|_{L^m \cap L^q}$ and $\big\||u_t(\tau,\cdot)|^p\big\|_{L^ q}$ will be handled as before to obtain (\ref{t6A1}). In order to control the integral with $\big\||u_t(\tau,\cdot)|^p\big\|_{\dot{H}^{s-2\delta, q}}$, we shall apply Corollary \ref{Corfractionalhomogeneous} for fractional powers with $s-2\delta \in \big(\frac{n}{q},p\big)$ and Corollary \ref{CorollaryEmbedding} with a suitable $s^* <\frac{n}{q}$. Hence, we have
$$ \big\||u_t(\tau,\cdot)|^p\big\|_{\dot{H}^{s-2\delta, q}} \lesssim \|u_t(\tau,\cdot)\|_{\dot{H}^{s-2\delta, q}}\|u_t(\tau,\cdot)\|^{p-1}_{L^\ity} \lesssim \|u_t(\tau,\cdot)\|_{\dot{H}^{s-2\delta, q}}\big(\|u_t(\tau,\cdot)\|_{\dot{H}^{s^*, q}}+ \|u_t(\tau,\cdot)\|_{\dot{H}^{s-2\delta, q}}\big)^{p-1}. $$
Applying the fractional Gagliardo-Nirenberg inequality leads to
$$ \|u_t(\tau,\cdot)\|_{\dot{H}^{s^*, q}} \lesssim \|u_t(\tau,\cdot)\|^{1-\theta}_{L^q}\,\,\big\||D|^{s-2\delta} u_t(\tau,\cdot)\big\|^{\theta}_{L^q} \lesssim (1+\tau)^{(2+[\frac{n}{2}])(1-\frac{\sigma}{2\delta})\frac{1}{r} -\frac{n}{2\delta}(1-\frac{1}{r}) -\frac{s^*}{2\delta}}\|u\|_{X(\tau)}, $$
where $\theta= \frac{s^*}{s-2\delta}$. Therefore, we obtain
\begin{align*}
&\big\||u_t(\tau,\cdot)|^p\big\|_{\dot{H}^{s-2\delta, q}} \lesssim (1+\tau)^{p\big((2+[\frac{n}{2}])(1-\frac{\sigma}{2\delta})\frac{1}{r} -\frac{n}{2\delta}(1-\frac{1}{r})\big)-\frac{s-2\delta}{2\delta}- (p-1)\frac{s^*}{2\delta}}\|u\|^p_{X(\tau)}\\
&\qquad \lesssim (1+\tau)^{p((2+[\frac{n}{2}])(1-\frac{\sigma}{2\delta})\frac{1}{r}-\frac{n}{2\delta}(\frac{1}{m}-\frac{1}{mp}))}\|u\|^p_{X(\tau)},
\end{align*}
if we choose $s^*= \frac{n}{q}- \e$, where $\e$ is a sufficiently small positive number. By the same arguments as above it follows
\begin{equation}
\big\||D|^{s-2\delta} u_t^{nl}(t,\cdot)\big\|_{L^q} \lesssim (1+t)^{1+(2+[\frac{n}{2}])(1-\frac{\sigma}{2\delta})\frac{1}{r} -\frac{n}{2\delta}(1-\frac{1}{r}) -\frac{s}{2\delta}} \|u\|^p_{X(t)}. \label{t6A3}
\end{equation}
Analogously, we also get
\begin{equation}
\big\||D|^s u^{nl}(t,\cdot)\big\|_{L^q} \lesssim (1+t)^{1+(1+[\frac{n}{2}])(1-\frac{\sigma}{2\delta})\frac{1}{r} -\frac{n}{2\delta}(1-\frac{1}{r}) -\frac{s}{2\delta}} \|u\|^p_{X(t)}. \label{t6A4}
\end{equation}
From (\ref{t6A1}) to (\ref{t6A4}) and the definition of the norm in $X(t)$ we may conclude immediately the inequality (\ref{pt4.31}).\medskip

\noindent Next, let us prove the inequality (\ref{pt4.4}). The new difficulty is to control the term $\big\||u_t(\tau,\cdot)|^p-|v_t(\tau,\cdot)|^p\big\|_{\dot{H}^{s-2\delta,q}}$. Then, repeating the proof of Theorem \ref{dl2.3} and using the analogous treatment as in the first step, we obtain
\begin{align*}
\big\||D|^{s-2\delta} \partial_t \big(Nu(t,\cdot)- Nv(t,\cdot)\big)\big\|_{L^q} &\lesssim  (1+t)^{1+(2+[\frac{n}{2}])(1-\frac{\sigma}{2\delta})\frac{1}{r} -\frac{n}{2\delta}(1-\frac{1}{r}) -\frac{s}{2\delta}} \|u-v\|_{X(t)}\big( \|u\|^{p-1}_{X(t)}+ \|v\|^{p-1}_{X(t)} \big), \\
\big\||D|^s \big(Nu(t,\cdot)- Nv(t,\cdot)\big)\big\|_{L^q} &\lesssim (1+t)^{1+(1+[\frac{n}{2}])(1-\frac{\sigma}{2\delta})\frac{1}{r} -\frac{n}{2\delta}(1-\frac{1}{r}) -\frac{s}{2\delta}} \|u-v\|_{X(t)}\big( \|u\|^{p-1}_{X(t)}+ \|v\|^{p-1}_{X(t)} \big).
\end{align*}
From the definition of the norm in $X(t)$ we may conclude immediately the inequality (\ref{pt4.4}). This completes the proof of Theorem \ref{dl2.5}.

\section{Concluding remarks and open problems} \label{Sec5}

\begin{nx}{(Time-dependent coefficients in the dissipation term)}
\fontshape{n}
\selectfont
A next challenge is to study  $L^1$ estimates for oscillating integrals and $L^{p}- L^{q}$ linear estimates away from the conjugate line as well to structurally damped $\sigma$-evolution models with time-dependent coefficients. These estimates are fundamental tools to prove global (in time) existence results to semi-linear models. Therefore, it is interesting to investigate the following Cauchy problem:
\begin{equation}
u_{tt}+ (-\Delta)^\sigma u+ b(t) (-\Delta)^\delta u_t = 0,\, u(0,x)= u_0(x),\, u_t(0,x)=u_1(x) \label{pt0.1}
\end{equation}
with $\sigma \ge 1$ and $\delta \in [0,\sigma]$. Here the coefficient $b=b(t)$ should satisfy some ``effectiveness assumptions'' as in \cite{Kainane}.
\end{nx}

\begin{nx}{(Blow-up results)}
\fontshape{n}
\selectfont
In this paper, we applied $(L^{m}\cap L^{q})- L^{q}$ and $L^{q}- L^{q}$ estimates for solution and its derivatives to (\ref{pt6.3}) to prove the global (in time) existence of small data Sobolev solutions to the semi-linear models (\ref{pt6.1}) and (\ref{pt6.2}) with $\delta \in (\frac{\sigma}{2},\sigma]$. It can be expected to find the critical exponents for each of the two nonlinearities. The ``shape" of these critical exponents can be found in \cite{DabbiccoEbert} by using the test function method, where the assumption for integers $\sigma$ and $\delta$ comes into play. In general, the main difficulty is to deal with fractional Laplacian operators $(-\Delta)^\sigma$ as well-known non-local operators. 
\end{nx}

\begin{nx}{(Gevrey smoothing)}
\fontshape{n}
\selectfont
We are interested in another qualitative property of solutions to (\ref{pt6.3}), the so-called Gevrey smoothing. It is reasonable to use our estimates with $L^2$ norms only. Moreover, we suppose for the Cauchy data $(u_0,u_1) \in H^\sigma \times L^2$. The study of regularity properties for the solutions allows to restrict our considerations to large frequencies in the extended phase space. Recalling the definition the Gevrey-Sobolev space regularity $\Gamma^{s,\rho}$ introduced in \cite{DaoReissig} we may conclude the following statement.\medskip

\noindent \textbf{Theorem 6.} \textit{Let us consider the Cauchy problem (\ref{pt6.3}) with $\delta\in (\frac{\sigma}{2},\sigma)$. The data are supposed to belong to the energy space, that is, $(u_0,u_1) \in H^\sigma \times L^2$. Then, there is a smoothing effect in the sense, that the solution belongs to the Gevrey-Sobolev space as follows:
\[ u(t,\cdot) \in \Gamma^{\frac{1}{2(\sigma-\delta)},\sigma}, \text{ and } |D|^\sigma u(t,\cdot),\,\, u_t(t,\cdot) \in \Gamma^{\frac{1}{2(\sigma-\delta)},0} \text{ for all } t>0. \]}

\begin{proof}
Let us turn to large values of $|\xi|$. For the sake of the asymptotic behavior of the characteristic roots in (\ref{pt3.3}) and (\ref{pt3.4}) we arrive at
$$ |\hat{K_0}|\lesssim e^{-c|\xi|^{2(\sigma-\delta)}t},\,\, |\hat{K_1}|\lesssim |\xi|^{-2\delta}e^{-c|\xi|^{2(\sigma-\delta)}t} \text{ and } |\partial_t \hat{K_0}|\lesssim |\xi|^{2(\sigma-\delta)}e^{-c|\xi|^{2(\sigma-\delta)}t},\,\, |\partial_t \hat{K_1}|\lesssim e^{-c|\xi|^{2(\sigma-\delta)}t}, $$
for some positive constants $c$. Therefore, using the representation of the solutions (\ref{pt3.2}) we derive the following estimates:
\begin{align*}
\int_{\R^n}\exp \big(2c|\xi|^{2(\sigma-\delta)}t \big) |\xi|^{2\sigma}|v(t,\xi)|^2 d\xi &\lesssim \int_{\R^n}|\xi|^{2\sigma} |v_0(\xi)|^2 d\xi + \int_{\R^n}|v_1(\xi)|^2 d\xi, \\
\int_{\R^n}\exp \big(2c|\xi|^{2(\sigma-\delta)}t \big) |v_t(t,\xi)|^2 d\xi &\lesssim \int_{\R^n}|\xi|^{2\sigma} |v_0(\xi)|^2 d\xi + \int_{\R^n}|v_1(\xi)|^2 d\xi.
\end{align*}
We may conclude immediately all the statements we wanted to prove.
\end{proof}
\end{nx}




\begin{thebibliography}{00}

\bibitem{DabbiccoEbert2014} M. D’Abbicco, M.R. Ebert, \textit{An application of $L^{p}-L^{q}$ decay estimates to the semilinear wave equation with parabolic-like structural damping}, Nonlinear Analysis, 99 (2014), 16-34.
\bibitem{DabbiccoEbert} M. D’Abbicco, M.R. Ebert, \textit{A new phenomenon in the critical exponent for structurally damped semi-linear evolution equations}, Nonlinear Analysis, 149 (2017), 1-40.
\bibitem{DuongKainaneReissig} Duong T. P., M. Kainane Mezadek, and M. Reissig, \textit{Global existence for semi-linear structurally damped $\sigma$-evolution models}, J. Math. Anal. Appl., 431 (2015), 569-596.
\bibitem{DabbiccoReissig} M. D'Abbicco, M. Reissig, \textit{Semilinear structural damped waves}, Math. Methods Appl. Sci., 37 (2014), 1570-1592.
\bibitem{DaoReissig} T.A. Dao, M. Reissig, \textit{An application of $L^1$ estimates for oscillating integrals to parabolic like semi-linear structurally damped $\sigma$-evolution models}, 33A4, submitted.
\bibitem{ReissigEbert} M. R. Ebert, M. Reissig, ``Methods for partial differential equations, qualitative properties of solutions, phase space analysis, semilinear models'',  Birkh\"{a}user, 2018.
\bibitem{Grafakos} L. Grafakos, ``Classical and modern Fourier analysis'', Prentice Hall, 2004.
\bibitem{GalaktionovMitidieri} V.A. Galaktionov, E.L. Mitidieri, and S.I. Pohozaev, ``Blow-up for higher-order prabolic, hyperbolic, dispersion and Schr\"{o}dinger equations'', in: Monogr. Res. Notes Math., Chapman and Hall/CRC, ISBN: 9781482251722, 2014.
\bibitem{Ozawa} H. Hajaiej, L. Molinet, T. Ozawa, and B. Wang,  ``Necessary and sufficient conditions for the fractional Gagliardo-Nirenberg inequalities and applications to Navier-Stokes and generalized boson equations, Harmonic analysis and nonlinear partial differential equations'', 159-175, RIMS Kokyuroku Bessatsu, B26, Res.Inst.Math.Sci. (RIMS), Kyoto, 2011.
\bibitem{Ikehata} R. Ikehata, \textit{Asymptotic profiles for wave equations with strong damping}, J. Differential Equations, 257 (2014) 2159-2177.
\bibitem{IkehataTodorova} R. Ikehata, G. Todorova, and B. Yordanov, \textit{Wave equations with strong damping in Hilbert spaces, Journal Differential Equations}, 254 (2013), 3352-3368.
\bibitem{Kainane} M. Kainane, ``Structural damped $\sigma$-evolution operators'', PhD thesis, TU Bergakademie Freiberg, Germany, 2013.
\bibitem{Marcinkiewicz} J. Marcinkiewicz, \textit{Sur les multiplicateurs des s\'{e}ries de Fourier}, Studia Math, 8 (1939), 78-91.
\bibitem{Miyachi} A. Miyachi, \textit{On some Fourier multipliers for $H^p(\R^n)$}, J. Fac. Sci. Univ. Tokyo IA, 27 (1980), 157-179.
\bibitem{Miyachi1980} A. Miyachi, \textit{On some estimates for the wave equation in $L^p$ and $H^p$}, J. Fac. Sci. Univ. Tokyo IA, 27 (1980), 331-354.
\bibitem{MitidieriPohozaev} E. Mitidieri, S.I. Pohozaev, \textit{Non-existence of weak solutions for some degenerate elliptic and parabolic problems on $\R^n$}, J. Evol. Equ., 1 (2001), 189-220.
\bibitem{Bui} Bui Tang Bao Ngoc, ``Semi-linear waves with time-pendent speed and dissipation'', PhD thesis, TU Bergakademie Freiberg, Germany, 2014.
\bibitem{NarazakiReissig} T. Narazaki, M. Reissig, \textit{$L^1$ estimates for oscillating integrals related to structural damped wave models}, in: M. Cicognani, F. Colombini, D. Del Santo (Eds.), Studies in Phase Space Analysis with Applications to PDEs, in: Progr. Nonlinear Differential Equations Appl., Birkh\"{a}user, 2013, 215-258.
\bibitem{Peral} J.C. Peral, \textit{$L^p$ estimates for the wave equation}, J. Funct. Anal., 36 (1980), 114-145.
\bibitem{Pizichillo} F. Pizichillo, ``Linear and non-linear damped wave equations'', Master thesis, 62pp., University of Bari, 2014.
\bibitem{Palmierithesis} A. Palmieri, M. Reissig, \textit{Semi-linear wave models with power non-linearity and scale-invariant time-dependent mass and dissipation, II}, Mathematische Nachrichten., (2018), 1-34, https://doi.org/10.1002/mana.201700144.
\bibitem{RunSic} T. Runst,  W. Sickel, ``Sobolev spaces of fractional order, Nemytskij operators, and nonlinear partial differential equations, De Gruyter series in nonlinear analysis and applications'', Walter de Gruyter $\&$ Co., Berlin, 1996.
\bibitem{FrancescoBruno} Cav. Francesco Fa\`{a} di Bruno, \textit{Note sur une nouvelle formule de calcul differentiel}, Quarterly J. Pure Appl. Math., 1 (1857), 359-360.
\bibitem{Simander} C. G. Simander, ``On Dirichlet boundary value problem", An $L^p$-Theory based on a generalization of G{\aa}rding's inequality, Lecture Notes in Mathematics, 268, Springer, Berlin, 1972.
\bibitem{Shibata} Y. Shibata, \textit{On the rate of decay of solutions to linear viscoelastic equation}, Math. Methods Appl. Sci., 23 (2000), 203-226.
\bibitem{SteinWeiss} E. Stein, G. Weiss, \textit{Fractional integrals on $n$-dimensional Euclidean space}, J. Math. Mech., 7 (1958), 503-514.
\bibitem{Weisz} F. Weisz, \textit{Marcinkiewicz multiplier theorem and the Sunouchi operator for Ciesielski–Fourier series}, Journal of Approximation Theory, 133 (2005), 195-220.

\end{thebibliography}


\noindent \textbf{Acknowledgments}\medskip

\noindent The PhD study of MSc. T.A. Dao is supported by Vietnamese Government's Scholarship.\medskip

\noindent\textbf{Appendix A}\medskip

\noindent \textit{A.1. Fractional Gagliardo-Nirenberg inequality}

\begin{md} \label{fractionalGagliardoNirenberg}
Let $1<p,p_0,p_1<\infty$, $\sigma >0$ and $s\in [0,\sigma)$. Then, it holds the following fractional Gagliardo-Nirenberg inequality for all $u\in L^{p_0} \cap \dot{H}^\sigma_{p_1}$:
$$ \|u\|_{\dot{H}^{s}_p} \lesssim \|u\|_{L^{p_0}}^{1-\theta}\|u\|_{\dot{H}^{\sigma}_{p_1}}^\theta, $$
where $\theta=\theta_{s,\sigma}(p,p_0,p_1)=\frac{\frac{1}{p_0}-\frac{1}{p}+\frac{s}{n}}{\frac{1}{p_0}-\frac{1}{p_1}+\frac{\sigma}{n}}$ and $\frac{s}{\sigma}\leq \theta\leq 1$ .
\end{md}
For the proof one can see \cite{Ozawa}.
\medskip

\noindent \textit{A.2. Fractional Leibniz rule}

\begin{md} \label{fractionalLeibniz}
Let us assume $s>0$ and $1\leq r \leq \infty, 1<p_1,p_2,q_1,q_2 \le \infty$ satisfying the relation \[ \frac{1}{r}=\frac{1}{p_1}+\frac{1}{p_2}=\frac{1}{q_1}+\frac{1}{q_2}.\]
Then, the following fractional Leibniz rule holds:
$$ \|\,|D|^s(u \,v)\|_{L^r}\lesssim \|\,|D|^s u\|_{L^{p_1}}\|v\|_{L^{p_2}}+\|u\|_{L^{q_1}}\|\,|D|^s v\|_{L^{q_2}} $$
for any $u\in \dot{H}^s_{p_1} \cap L^{q_1}$ and $v\in \dot{H}^s_{q_2} \cap L^{p_2}$.
\end{md}
These results can be found in \cite{Grafakos}.
\medskip

\noindent \textit{A.3. Fractional chain rule}

\begin{md} \label{Propfractionalchainrulegeneral}
Let us choose $s>0$, $p>\lceil s \rceil$
 and $1<r,r_1,r_2<\infty$ satisfying $\frac{1}{r}=\frac{p-1}{r_1}+\frac{1}{r_2}$. Let us denote by $F(u)$ one of the functions $|u|^p, \pm |u|^{p-1}u$. Then, it holds the following fractional chain rule:
$$ \|\,|D|^{s} F(u)\|_{L^r}\lesssim \|u\|_{L^{r_1}}^{p-1}\|\,|D|^{s} u\|_{L^{r_2}} $$
for any $u\in  L^{r_1} \cap \dot{H}^{s}_{r_2}$.
 \end{md}
The proof can be found in \cite{Palmierithesis}.
\medskip

\noindent \textit{A.4. Fractional powers}

\begin{md} \label{PropSickelfractional}
Let $p>1$, $1< r <\infty$ and $u \in H^{s}_r$, where $s \in \big(\frac{n}{r},p\big)$.
Let us denote by $F(u)$ one of the functions $|u|^p,\, \pm |u|^{p-1}u$ with $p>1$. Then, the following estimate holds$:$
$$\Vert F(u)\Vert_{H^{s}_r}\lesssim \|u\|_{H^{s}_r}\|u\|_{L^\infty}^{p-1}.$$
\end{md}

\begin{hq} \label{Corfractionalhomogeneous}
Under the assumptions of Proposition \ref{PropSickelfractional} it holds: $\| F(u)\|_{\dot{H}^{s}_r}\lesssim \| u\|_{\dot{H}^{s}_r}\|u\|_{L^\infty}^{p-1}.$
\end{hq}
The proof can be found in \cite{DuongKainaneReissig}.
\medskip

\noindent \textit{A.5. A fractional Sobolev embedding}

\bmd \label{FracSobolevEmbedding}
Let $n \ge 1$, $0< s< n$, $1< q\le r< \ity$, $\alpha< \frac{n}{q^{'}}$ where $q^{'}$ denotes conjugate number of $q$, and $\gamma > -\frac{n}{r}$, $\alpha \ge \gamma$ satisfying $\frac{1}{r}= \frac{1}{q}+ \frac{\alpha- \gamma- s}{n}$. Then, it holds:
$$ \big\||x|^\gamma |D|^{-s}u\big\|_{L^r} \lesssim \big\||x|^\alpha u\big\|_{L^q}, \text{ that is }, \big\||x|^\gamma u\big\|_{L^r} \lesssim \big\||x|^\alpha\, |D|^s u\big\|_{L^q} $$
for any $u \in \dot{H}^{s,q}_\alpha$, where $\dot{H}^{s,q}_\alpha= \{u \,:\, |D|^s u \in L^q(\R^n,\, |x|^{\alpha q})\}$ is the weighted homogeneous Sobolev space of potential type with the norm $\|u\|_{\dot{H}^{s,q}_\alpha}= \big\||x|^\alpha\, |D|^s u \big\|_{L^q}$.
\emd
The proof can be found in \cite{SteinWeiss}.

\bhq \label{CorollaryEmbedding}
Let $1< q< \ity$ and $0< s_1< \frac{n}{q}< s_2$. Then, for any function $u \in \dot{H}^{s_1,q} \cap \dot{H}^{s_2,q}$ we have
\[ \|u\|_{L^\ity} \lesssim \|u\|_{\dot{H}^{s_1,q}}+ \|u\|_{\dot{H}^{s_2,q}}. \]
\ehq
For the proof one can see \cite{DaoReissig}.\medskip

\noindent \textit{A.6. A variant of Mikhlin- H\"{o}mander multiplier theorem}

\begin{md} \label{PropositionMultiplier} 
Let $q\in (1,\ity)$, $k=[\frac{n}{2}]+1$ and $b\ge 0$. Suppose that $m\in C^k(\R^n)$ satisfies $m(\xi)=0$ if $|\xi|\le 1$ and
$$ \big|\partial_\xi^\alpha m(\xi)\big| \le C|\xi|^{-nb |\frac{1}{q}-\frac{1}{2}|}\big(A|\xi|^{b-1}\big)^{|\alpha|} $$
for all $|\alpha|\le k$, $|\xi|\ge 1$ and with some constants $A\ge 1$. Then, the operator $T_m=F^{-1}\big(m(\xi)\big) \ast_{(x)}$, defined by the action \[ T_m f(x):= F^{-1}_{\xi \rightarrow x}\big(m(\xi) F_{y \rightarrow \xi}\big(f(y)\big)\big),\] is continuously bounded from $L^q$ into itself and satisfies the following estimate:
\[ \|T_m f(\cdot)\|_{L^q} \le C A^{n  |\frac{1}{q}-\frac{1}{2}|}\|f\|_{L^q}. \]
\end{md}
The proof of this lemma can be found in \cite{DabbiccoEbert} (Theorem $10$) and \cite{Miyachi} (Theorem $1$).\medskip

\noindent \textit{A.7. Modified Bessel functions}

\begin{md} \label{FourierModifiedBesselfunctions}
Let $f \in L^p(\R^n)$, $p\in [1,2]$, be a radial function. Then, the Fourier transform $F(f)$ is also a radial function and it satisfies
$$ F_n(\xi):= F(f)(\xi)= c \int_0^\ity g(r) r^{n-1} \tilde{J}_{\frac{n}{2}-1}(r|\xi|)dr,\,\,\, g(|x|):= f(x), $$
where $\tilde{J}_\mu(s):=\frac{J_\mu(s)}{s^\mu}$ is called the modified Bessel function with the Bessel function $J_\mu(s)$ and a non-negative integer $\mu$.
\end{md}

\begin{md} \label{PropertiesModifiedBesselfunctions}
The the following properties of the modified Bessel function hold:
\begin{enumerate}
\item $sd_s\tilde{J}_\mu(s)= \tilde{J}_{\mu-1}(s)-2\mu \tilde{J}_\mu(s)$,
\item $d_s\tilde{J}_\mu(s)= -s\tilde{J}_{\mu+1}(s)$,
\item $\tilde{J}_{-\frac{1}{2}}(s)= \sqrt{\frac{2}{\pi}}\cos s$ and $\tilde{J}_{\frac{1}{2}}(s)= \sqrt{\frac{2}{\pi}} \frac{\sin s}{s}$,
\item $|\tilde{J}_\mu(s)| \le Ce^{\pi|Im\mu|} \text{ if } s \le 1, $ and $\tilde{J}_\mu(s)= Cs^{-\frac{1}{2}}\cos \big( s-\frac{\mu}{2}\pi- \frac{\pi}{4} \big) +\mathcal{O}(|s|^{-\frac{3}{2}}) \text{ if } |s|\ge 1$,
\item $\tilde{J}_{\mu+1}(r|x|)= -\frac{1}{r|x|^2}\partial_r \tilde{J}_\mu(r|x|)$, $r \ne 0$, $x \ne 0$.
\end{enumerate}
\end{md}
\medskip

\noindent \textit{A.8. Fa\`{a} di Bruno's formula}

\begin{md} \label{FadiBruno'sformula1}
Let $h\big(g(x)\big)= (h\circ g)(x)$ with $x\in \R$. Then, we have
$$ \frac{d^n}{dx^n}h\big(g(x)\big)= \sum \frac{n!}{m_1! 1!^{m_1} m_2! 2!^{m_2}\cdots m_n! n!^{m_n}}h^{(m_1+m_2+\cdots+m_n)}\big(g(x)\big) \prod_{j=1}^n \big( g^{(j)}(x) \big)^{m_j}, $$
where the sum is taken over all $n$- tuples of non-negative integers $(m_1,m_2,\cdots,m_n)$ satisfying the constraint of the Diophantine equation: $ 1\cdot m_1+ 2\cdot m_2+\cdots+ n\cdot m_n =n.$
\end{md}
For the proof one can see \cite{FrancescoBruno}.
\medskip

\noindent \textit{A.9. A useful lemma}

\bbd \label{LemmaDerivative}
The following formula of derivative of composed function holds for any multi-index $\alpha$:
$$ \partial_\xi^\alpha h\big(f(\xi)\big)= \sum_{k=1}^{|\alpha|}h^{(k)} \big(f(\xi)\big)\Big(\sum_{\substack{\gamma_1+\cdots+\gamma_k \le \alpha\\ |\gamma_1|+\cdots+|\gamma_k|= |\alpha|,\, |\gamma_i|\ge 1}}\big(\partial_\xi^{\gamma_1} f(\xi)\big) \cdots \big(\partial_\xi^{\gamma_k} f(\xi)\big)\Big), $$
where $h=h(s)$ and $h^{(k)}(s)=\frac{d^k h(s)}{ds^k}$.
\ebd
The result can be found in \cite{Simander} at the page $202$.\medskip

\bbd \label{LemmaIntegral}
Let $\alpha, \beta \in \R.$ Then:
$$ I(t):=\int_0^t (1+t-\tau)^{-\alpha}(1+\tau)^{-\beta}d\tau \lesssim
\begin{cases}
(1+t)^{-\min\{\alpha, \beta\}} \hspace{2cm} \text{ if } \max\{\alpha, \beta\}>1,&\\
(1+t)^{-\min\{\alpha, \beta\}}\log(2+t) \hspace{0.45cm} \text{ if } \max\{\alpha, \beta\}=1,&\\
(1+t)^{1-\alpha-\beta} \hspace{2.5cm} \text{ if } \max\{\alpha, \beta\}<1.&\\
\end{cases} $$
\ebd
For the proof one can see \cite{DaoReissig}.


\end{document}